 \documentclass[final,1p,times]{elsarticle}

\usepackage{amssymb}

\usepackage{amsmath}
\usepackage{graphicx}
\usepackage[colorlinks,citecolor=red,urlcolor=blue,bookmarks=false,hypertexnames=true]{hyperref}

\usepackage{bbm,bm,dsfont,subfigure,multirow,tikz}
\usepackage{amsmath}
\usepackage{bm}
\usepackage{multicol}
\usepackage{float}
\usepackage{color}
\usepackage{mathrsfs}
\newtheorem{theorem}{Theorem}
\newtheorem{lemma}{Lemma}
\newtheorem{remark}{Remark}

\newtheorem{example}{Example}
\usepackage{algorithmicx,algorithm}
\usepackage{algpseudocode}
\newenvironment{proof}{{\noindent\it Proof}\quad}{\hfill $\square$\par}  
\numberwithin{equation}{section}
\numberwithin{remark}{section}
\numberwithin{theorem}{section}
\numberwithin{lemma}{section}
\numberwithin{corollary}{section}
\numberwithin{assumption}{section}
\numberwithin{figure}{section}
\numberwithin{table}{section}
\numberwithin{example}{section}

\def\tc{\text{curl}}

\def\NN{\mathbb{N}}

\def\RR{\mathbb{R}}

\def\NN{\mathbb{N}}
\def\curl{\text{curl}}
\def\rot{\text{rot}}

\def\tr{\mathsf{T}}

\DeclareMathOperator{\diag}{diag}

\def\IIone{\rm\uppercase\expandafter{\romannumeral1}}
\def\IItwo{\rm\uppercase\expandafter{\romannumeral2}}
\def\IIthree{\rm\uppercase\expandafter{\romannumeral3}}

\journal{Computers and Mathematics with Applications}

\begin{document}

\begin{frontmatter}



\title{Fast Maxwell Solvers Based on Exact Discrete Eigen-Decompositions I. Two-Dimensional Case}


%

 \author[label1]{Lixiu Wang}
 \ead{lxwang@ustb.edu.cn}
 
 \author[label2]{Lueling Jia\corref{cor1}}
 \ead{lljia@sdnu.edu.cn}
 
 \author[label3,label4]{Zijian Cao}
 \ead{2830769384@qq.com}
 
 \author[label4]{Huiyuan Li}
\ead{huiyuan@iscas.ac.cn}

\author[label5]{Zhimin Zhang}
\ead{ag7761@wayne.edu}

\cortext[cor1]{Corresponding author}

\address[label1]{School of Mathematics and Physics, University of Science and Technology Beijing, Beijing 100083, China.}

\address[label2]{School of Mathematics and Statistics, Shandong Normal University, Jinan 250358, China. }

\address[label3]{School of Computer Science and Technology, University of Chinese Academy of Sciences, Beijing 101408, China.}

\address[label4]{ State Key Laboratory of Computer Science/Laboratory of Parallel Computing, Institute of Software, Chinese Academy of Sciences, Beijing 100190, China.}

\address[label5]{Department of Mathematics, Wayne State University, Detroit, MI 48202, USA.}

\begin{abstract}
In this paper, we  propose fast solvers for Maxwell's equations in {\color{black}rectangular domains}.
		We first discretize the simplified  Maxwell's eigenvalue  problems  by employing the lowest-order rectangular N\'ed\'elec elements and derive the discrete eigen-solutions explicitly,  providing  a Hodge-Helmholtz decomposition framework at the discrete level.
		Based on exact eigen-decompositions, we further  design fast solvers for various Maxwell's source problems, guaranteeing either the divergence-free constraint or the Gauss's law at the discrete level.  With the help of fast sine/cosine transforms, the computational time  grows asymptotically as $\mathcal{O}(n^2\log n)$ with $n$ being the number of grids in each direction.  Our fast Maxwell solvers
		outperform other existing Maxwell solvers in  the literature and 
		 fully rival fast scalar Poisson/Helmholtz solvers based  on trigonometric transforms  in either efficiency, robustness, or storage complexity.
		It is also utilized to perform an efficient pre-conditioning for solving Maxwell's source problems with variable coefficients.
		Finally, numerical experiments are carried out to illustrate the effectiveness and efficiency of the proposed fast solver.

\end{abstract}




\end{frontmatter}



\section{Introduction}
	Electromagnetic phenomena are described by Maxwell's equations \cite{Orfanidis2016Electromagnetic},
	 originally presented using first-order curl operators. Many of the Maxwell's equations can be transformed into the curl-curl formulation.
	Consider the  model problem of the  Maxwell's curl-curl  equations,
	\begin{subequations}
	\begin{align}
	\label{Maxwell_source}
			&{\color{black}\curl\ \rot\  \bm{u}} +\alpha \bm u= \bm{f} \quad {\rm in}\;\Omega,
			\\
        	\label{divergence}			
			&\nabla\cdot \bm{u} = 0\quad {\rm in}\; \Omega,
	\end{align}
	\end{subequations}
	subject to  one of the   following essential/natural   boundary conditions
		\begin{subequations}
	\begin{align}
           \label{essential_BCs}			
			&\bm u \times \bm n=0  \quad\text {on } \partial\Omega,
			\\
\label{natural_BCs}
		&{\color{black} \rot\ }\bm u =0\; {\text{and}} \;\bm u\cdot \bm n=0\quad\text{on } \partial\Omega,	
	\end{align}
	\end{subequations}	
	where  $\bm u$  represents an electric field, $\bm f$ is  the current density,  $\alpha\in\RR$ is  a given constant,  and $\bm n$ is the unit outward normal of the boundary $\partial \Omega$.  	  
	
	
	The curl-curl equations \eqref{Maxwell_source}-\eqref{divergence}  are widely encountered in  
	time-dependent electromagnetic fields \cite{Li2013Time}, time-harmonic fields \cite{Sun2016Finite, Liuqiang2009OPTIMAL}, cavity eigenvalue problems \cite{Monk2003},  magnetohydrodynamic problems \cite{LI2021109980A}, and so on.  	
	Therefore,  the development of fast solvers along with  appropriate numerical approximation schemes is of particular importance. 

	N\'ed\'elec originally   proposed two families of curl-conforming edge elements \cite{Lec1980A,Lec1986A}  in 1980s, which exactly satisfy the intrinsic requirement  of  tangential continuity of an electric field.  Kikuchi  \cite{Kikuchi1987Mixed}  introduced  a mixed formulation by  incorporating a Lagrange multiplier that  maintains the  divergence-free constraint  in the weak sense. Their efforts laid a state-of-art foundation for Maxwell's curl-curl equations,  guaranteeing  {\color{black}the well-posedness of discrete problems \cite{Sun2016Finite, Kikuchi1989Mixed} }and  maintaining the divergence-free constraint at  a certain discrete level.  {\color{black}However, the resulting linear algebraic system is larger and of 
	 saddle-point type, making it more difficult to solve than the elliptic problem}. Moreover,  it may induce instabilities and  spurious solutions,  particularly for eigenvalue problems, {\color{black}as the} degrees  of freedom increase.

 To solve the saddle-point problems effectively,  some methods have been proposed, such as  null-space methods \cite{Rees2018A, Benzi2005Numerical} for transforming an indefinite system into a symmetric positive definite one of smaller dimension, coupled direct solvers  \cite{Benzi2005Numerical} based on triangular factorizations, and Krylov subspace methods \cite{vandervorst_2003Iterative}.  Besides,   preconditioning techniques also have been widely studied, such as block triangular preconditioners \cite{Liu2013New,Elman2002Performance},  inexact constraint preconditioners \cite{Bai2005On}, shift-splitting preconditioners \cite{Benzi2004A,Salkuyeh2015A}, primal-based penalty preconditioners \cite{Dohrmann2006A} and the references therein.  
	Moreover, to avoid the saddle-point problem, a $\delta$-regularization method was  developed in \cite{Duan2012A} to completely ignore the divergence-free constraint. 
  
  {\color{black}  There are some work  on fast solvers for Maxwell's curl-curl  equations. The discrete H(curl) systems could be efficiently solved by geometric multigrids  \cite{Multigrid2000Arnold, Multigrid1998Hiptmair,Xu2009Optimal}, as well as the nodal auxiliary space preconditioner  \cite{Hiptmair2007Nodal, Hu2021Convergence, Nodal2024Li}, i.e., the popular Hiptmair–Xu  preconditioner. Besides, the domain decomposition methods
  are frequently employed for solving problems in  H(curl), see \cite{Hu2003A, 2015PhDT220C, AFETI2001Andrea, Dohrmann2012An, MR3242973, Calvo2015ABA} and the reference therein. 
Although a fast algorithm  for three-dimensional photonic crystals was found in \cite{Lyu2021FAME},
		which combines the null-space free method with fast Fourier transform-based matrix-vector multiplications,
		little work has been reported on fast direct solvers for Maxwell's curl-curl  equations  up to the present.}

		The primary purpose  of this paper is to  design some fast  solvers  for the Maxwell's equations  \eqref{Maxwell_source}-\eqref{divergence} with  either \eqref{essential_BCs} or \eqref{natural_BCs}  for boundary conditions on a rectangle in  two dimensions. Distinct to existing work, the  Maxwell's curl-curl  equations are discretized 
		by using  the lowest-order edge elements with a naive variational form instead of the mixed form, and the resulting algebraic linear system  
		is then solved 
		 in the Hodge-Helmholtz decomposition framework   which is constructed through exact discrete eigen-decompositions.
		 Indeed, as proved in Section \ref{exactsolution} and Section \ref{Eigen-natural}, the discrete eigen-decompositions of the eigenvalue problem of the Maxwell's equation, 
		 \begin{align}
	\label{Maxwell_eigen}
			&{\color{black}\curl\ \rot\  \bm{u}} = \lambda \bm{u} \quad {\rm in}\;\Omega,
	\end{align}
	with 	boundary conditions either \eqref{essential_BCs} or \eqref{natural_BCs},  establishes two types of eigenfunctions constituting   the   entire approximation space:
	one  is divergence-free
		and corresponds to non-zero eigenvalues, while the other satisfies the curl-free constraint and associates with zero eigenvalues. 
	It is then not astonishing that the numerical solution of the Maxwell's source problem is simply a linear expansion of  all divergence-free  discrete eigenfunctions,
	whose coefficients can be efficiently evaluated  via fast subroutines such as FFTW \cite{fftw2020} for the discrete sine/cosine transforms (DSTs/DCTs) in $\mathcal{O}(n^2\log n)$ operations. In this way, fast Maxwell solvers are then designed just like fast solvers for scalar Poisson/Helmholtz equations. 
	Moreover, the divergence-free constraint is naturally guaranteed  at the discrete level, and instabilities and spurious solutions are then avoided. The obvious superiority of the fast Maxwell solvers on computational time and storage requirements over other existing direct solvers are demonstrated by various numerical experiments.

	We also aim at the fast solvers  guaranteeing   the Gauss's law at the discrete level for general Maxwell's equations,
	 \begin{subequations}
	\begin{align}\label{Maxwell_source_div}
			&{\color{black}\curl\ \rot\  \bm{u}} +\alpha \bm u +\nabla p= \bm{f}\quad \text{ in }\;\Omega,\\
			\label{Maxwell_Gauss_law}
			&\nabla\cdot \bm{u} = \rho\quad \text{ in }\; \Omega,\\
			\label{Maxwell_source_boundary}
			&\bm u \times \bm n=0, \;p=0 \quad \text{ on }\;\partial\Omega,
	\end{align}
	\end{subequations}
	where  $p$ represents the multiplier/pseudo-pressure and $\rho$  is the charge density.  Additionally, 
		an extended study on our fast solvers  as an efficient pre-conditioning for solving the Maxwell's source problem with variable coefficients
		is also carried out {\color{black}for $\alpha \ge  0$ and $ \beta \ge \beta_0$ for some $\beta_0 > 0$}.

	%
	%
	%

In summary, the novelty of this work resides in three aspects:
	\begin{itemize}
	 \item[(1)] {\color{black} By exploring the eigenvalue decomposition of the discrete Maxwell’s eigenvalue problem, we find a fast  way to compute the specific form of the Hodge-Helmholtz decomposition for electromagnetic fields on rectangular domains at the discrete level.}
	 
	 \item[(2)] We design fast solvers for various Maxwell's source problems  based on exact eigen-decompositions. The computational time  grows asymptotically as $\mathcal{O}(n^2\log n)$,  making our fast Maxwell solvers as efficient as the fast Poisson/Helmholtz solvers based on the trigonometric transforms.
	 
	\item[(3)] Our fast solver, based on the Hodge-Helmholtz decomposition, rigorously adheres to divergence-free constraint at the discrete level, thereby completely avoiding spurious solutions and numerical instability.
\end{itemize}

	%

	The remainder of this paper is organized as follows. In Section \ref{preli}, we formulate the approximation scheme for the curl-curl equations \eqref{Maxwell_source}-\eqref{divergence} and \eqref{essential_BCs}  in matrix form by specifying the basis functions of the lowest-order edge element space and the bilinear element space. Additionally, we provide some notations and preliminaries.
	In Section \ref{exactsolution}, we analyze  the exact  discrete eigen-decompositions of the Maxwell's eigenvalue problem \eqref{Maxwell_eigen}, and
	a fast solver for the discrete Maxwell's source problem is then designed in Section \ref{fastsolver}. 
	The algorithm is further extended to 
	the Maxwell's eigenvalue problem and source problems with  natural boundary conditions  in Section \ref{natural_bc}.  
Its extension  to general Maxwell's equations \eqref{Maxwell_source_div}-\eqref{Maxwell_source_boundary} is discussed  in Section \ref{solver_pressure}.
 Illustrative numerical results are reported in Section \ref{example} to confirm the computational efficiency.  Finally,  a conclusion remark is given in Section \ref{conclusion}.
	
	\section{Preliminaries}\label{preli}
	\subsection{Notations}
	Let ${\bm u}=(u_1, u_2)^{\tr}$ and ${\bm w}=(w_1, w_2)^{\tr}$ where the superscript $\tr$ denotes the transpose.
	Then ${\bm u} \times {\bm w} = u_1 w_2 - u_2 w_1$ and ${\color{black}\rot\ {\bm u}} = \partial u_2/\partial x - \partial u_1/\partial y$.
	For a scalar function $v$,   ${\color{black}\curl\ v} = (\partial v/\partial y, - \partial v/\partial x)^\tr$.
	Denote by $\NN$ and $\NN_0$ the set of positive integers and non-negative integers, respectively.
	For any $M,N\in\NN_0$, 
	we let $Q_{M,N}$ be the space of polynomials  whose separate degrees are no greater than $M$ in the $x$-direction and $N$ in the $y$-direction.
	
	Let ${\bm e}_k^n,\, 1\le k \le n$, be the $k$-th column of
	the identity matrix $I_{n}\in\RR^{n\times n}$ and simply denote $\bm e_0^n=\bm 0.$
	We further introduce the vectorization operation of a matrix \cite[Section 1.3.7]{Matrix2013}.
	Assume that $X=(\bm x_1,\bm x_2,\cdots,\bm x_p)^{\tr}\in \RR^{p\times q}$ is a matrix  such that $\bm x_i^\tr$ is the $i$-th row of $X$ for $ i=1,2,\cdots,p$.
	The vectorization of $X$,  denoted as $\overrightarrow{\text{vec}}(X)$, is a $pq$-by-1 vector obtained by stacking each row of $X$ and then transposing: $$\overrightarrow{\text{vec}}(X)=(\bm{x}_1^\tr, \bm{x}_2^\tr, \cdots,\bm{x}_p^\tr)^\tr.$$
	For ${A}\in \RR^{r\times p}$ and  ${B}\in\RR^{q\times s},$ it holds that
	\[\overrightarrow{\text{vec}}({A}X{B})=({A}\otimes {B}^\tr)\overrightarrow{\text{vec}}(X),\]
	where $\otimes$ represents the Kronecker product.
	
	
	\subsection{Approximation spaces and numerical schemes}
	\label{ASNS}
	For simplicity, we assume that $\Omega=[0,1]^2$ and let 
	$$\mathcal{T}_h=\left\{K_{i,j}:= [x_{i-1},x_{i}]\times[y_{j-1},y_{j}]: \ 1\le i,j\le n \right\}$$ be a $n\times n$ uniform mesh of $\Omega$ with $x_{\nu}=y_{\nu}=\nu h$ and $h=\frac{1}{n}$.
	Then the lowest-order edge element space is given by
	\begin{align}\label{U_h}
		&  U_h=\left\{\bm{w}_h\in H(\text{curl};\Omega):\ \bm w_h|_{K_{i,j}}\in Q_{0,1}\times Q_{1,0},\ \forall K_{i,j}\in\mathcal{T}_h\right\},
	\end{align}	
	where the  local N\'ed\'elec edge  basis functions are explicitly  presented in
	\cite{Li2013Time}. 
	Meanwhile, the bilinear element space is defined by
	\begin{align}\label{S_h}
		&  S_h=\left\{{q}_h\in H^1(\Omega):\  q_h|_{K_{i,j}}\in Q_{1,1},\; \forall K_{i,j}\in\mathcal{T}_h\right\}.
	\end{align}
Their subspaces with essential boundary conditions
 are given by 
	\begin{align}\label{dis_space_bd}
			&U^0_h= U_h \cap H_0(\tc;\Omega),
			&&S^0_h=S_h\cap H_0^1(\Omega),
	\end{align}
 {\color{black} where $H_0(\tc;\Omega):=\{\ \bm u\in H(\tc;\Omega): \ \bm u\times \bm n=0  \  \text{on}\ \partial\Omega\}.$}
   Then, the numerical scheme in the {\color{black}straightforward mixed} variational form for \eqref{Maxwell_source}-\eqref{essential_BCs}  
	  reads: to find $\bm{u}_h\in U_h^0$ such that 
	  \begin{subequations}
	\begin{align}\label{mix_source_scheme_a}
			&\left({\color{black} \rot\ } \bm{u}_h, {\color{black} \rot\ } \bm{w}_h\right) +\alpha(\bm u_h,\bm w_h)= (\bm{f},\bm{w}_h),\quad &\forall \bm{w}_h\in U_h^0,\\
			\label{mix_source_scheme_b}
			&\left(\bm{u}_h, \nabla q_h\right)=0,\quad &\forall q_h\in S_h^0.
	\end{align}
	\end{subequations}

 Let us write  $\bm u_h$  as
	\begin{align}\label{uh_series} 
		\bm u_h = \sum\limits_{j=1}^{n} \sum\limits_{i=1}^{n-1}  \widehat{u}_{i,j}\bm N^{(1)}_{i,j}(x,y) + 
		\sum\limits_{j=1}^{n-1} \sum\limits_{i=1}^{n}   \widehat{v}_{i,j}\bm N^{(2)}_{i,j}(x,y)\in U^0_h,
	\end{align}
	where   $\bm N^{(1)}_{i,j}(x,y)$ and $\bm N^{(2)}_{i,j}(x,y)$ are 
 the global basis functions  of $U_h^0$  corresponding to the (normalized) degrees of freedom  along the edge $\{x_i\}\times{\color{black}[} y_{j-1},y_j {\color{black}]}$ 
 and  the edge ${\color{black}[} x_{i-1},x_i {\color{black}]} \times\{y_j\}$  \cite{Li2013Time}, respectively.  {\color{black} The basis functions satisfy the tangent continuity, i.e.,
 \begin{align*}
    &\left(\bm N_{i,j}^{(1)}(x_m,y)\cdot \bm \tau\right)\Big|_{l_n}= \frac{1}{h}\delta_{i,m}\delta_{j,n},\ \text{for}\ l_n= [y_{n-1},y_n],\\
     &\left(\bm N_{i,j}^{(2)}(x,y_n)\cdot \bm \tau\right)\Big|_{l_m}= \frac{1}{h}\delta_{i,m}\delta_{j,n},\ \text{for}\ l_m= [x_{m-1},x_m],
\end{align*}
 where $\bm \tau$ is the anti-clock unit tangential vector.}

	Then we derive the equivalent linear algebraic system for the  approximation scheme \eqref{mix_source_scheme_a}-\eqref{mix_source_scheme_b}:
	\begin{subequations}
	\begin{align}\label{Kron-eq-01-sou}
		&\begin{cases}
			UA-B^\tr VB^\tr +\dfrac{ h^2}{6} \alpha UD= h^2 F,\\[0.5em]
			-BUB+AV +\dfrac{ h^2}{6}  \alpha DV=  h^2 G,
		\end{cases}
		\\
		\label{Kron-eq-02-sou}
		&BUD+DVB^\tr = \bm 0,
	\end{align}
	\end{subequations}
	where {\color{black}$U = \big( \widehat{u}_{i,j}\big)_{1\le i\le  n-1\atop 1\le j \le n}^\tr$,
	$V = \big( \widehat{v}_{i,j}\big)_{1\le i\le  n\atop 1\le j \le n-1}^\tr$,
	$F=\big( (\bm f,\bm N^{(1)}_{i,j}) \big)_{1\le i\le  n-1\atop 1\le j \le n}^\tr$,
	and $G=\big( (\bm f,\bm N^{(2)}_{i,j}) \big)_{1\le i\le  n\atop 1\le j \le n-1}^\tr$,}
	and  $A,D\in\RR^{(n-1)\times(n-1)}$ and $B\in\RR^{(n-1)\times n}$ are defined by 
	{\addtolength{\arraycolsep}{-0.1em}
		\begin{align}\label{importantmatrix}
			\begin{split}
				&A=\begin{pmatrix}
					2 &-1 &&\\
					-1 &2 &\ddots&\\
					&\ddots&\ddots&-1\\
					&&-1&2
				\end{pmatrix},\ 
				D=\begin{pmatrix}
					4 &1 &&\\
					1 &4 &\ddots&\\
					&\ddots&\ddots&1\\
					&&1&4
				\end{pmatrix},\ 
				B=\begin{pmatrix}
					-1 &1 &&&\\
					&-1&1&&\\
					& &\ddots&\ddots&\\
					&&&-1&1
				\end{pmatrix}. 
			\end{split}
		\end{align}
	}
 {\color{black} In fact, the matrices $A, B$ and $B^\tr$ correspond to the discrete one-dimensional Laplace operator,  divergence operator, and gradient operator, respectively. Moreover, it holds that $BB^\tr=A$.}

%

	\subsection{Discrete sine and cosine transforms}
	
	For $n\in\NN$, define the discrete sine transform (DST)  matrix 
	\begin{align}\label{matrixS}
		S=
		\left( \bm s_1, \bm s_2, \dots, \bm s_{n-1} \right),
  \end{align}
{\color{black}  where $\bm{s}_j =\Big(
		\sin \frac{j \pi}{n},    \sin \frac{2 j \pi}{n}, \dots,  \sin \frac{(n-1) j \pi}{n} 
		\Big)^\tr,\  0\le j\le n-1.$}
	
	For any $\bm x\in \RR^{n-1}$, the matrix-vector multiplication $S\bm x$ represents the DST of  the  first type (DST-\IIone) of $\bm x$.
	One readily checks that $S$ is symmetric {\color{black} and satisfies the following property  \cite{BRITANAK200773}}
	\begin{align}\label{S_inv}
		S^{-1} = \frac{2}{n}S.
	\end{align}
	Moreover, the diagonalizations of $A$ and $D$ defined in \eqref{importantmatrix} 
	are
	\begin{align}\label{diag_AD}
		S^{-1}A S=SA S^{-1}= \Lambda^2,\quad
		S^{-1}D S=SD S^{-1}= \Sigma,
	\end{align}	
	respectively, where	
	\begin{align}
		&\Lambda=\diag\left( \tau_1,\tau_2,\cdots,\tau_{n-1}\right), 
		\quad {\color{black}\tau_i = - \sqrt{ 2\left(1-\cos\frac{i\pi}{n}\right)}  =-2\sin \frac{i\pi}{2n}},
  \label{Lambda}
		\\
		&\Sigma = \diag\left( \sigma_1,\sigma_2,\cdots,\sigma_{n-1}\right),
		\quad\sigma_i = 2\Big( 2+\cos\frac{i\pi}{n}\Big).
		\label{Sigma}
	\end{align}	
	For $n\in\NN$, define the discrete cosine transform (DCT)  matrix 
	\begin{align}\label{matrixC2}
		\begin{split}
			&C= \left(  \bm c_0,\bm c_1,\cdots, \bm c_{n-1}\right),
			\end{split}
	\end{align}
		{\color{black} where $\bm c_j =\nu_j^n \Big(
			\cos\frac{j \pi}{2n},  \cos\frac{3 j \pi}{2n}, \cdots , \cos\frac{(2n-1)j \pi}{2n}
			\Big)^\tr,\ 0\le j\le n$   
	and 
	\begin{align}\label{nu_j}
		\nu_j^n=1 \text{ if } 1\le j\le n-1,\quad 
		\nu_j^n= \frac{1}{\sqrt{2}} \text{ if }  j=0 \text{ or }j=n.
	\end{align}}
	For any $\bm x\in \RR^{n\times 1}$, the matrix-vector multiplications $ C^\tr\bm x$ and $ C \bm x$ imply 
	the DCT of the second type (DCT-\IItwo) and of the third type (DCT-\IIthree) of $\bm x$, respectively.
	In addition, we have \cite{BRITANAK200773}
	\begin{align}\label{C_inv}
		C^{-1} = \frac{2}{n} C^{\tr}.
	\end{align}
	
	We emphasize that all DSTs and DCTs  can be implemented by fast sine/cosine transforms (FSTs/FCTs) in $\mathcal{O}(n \log n)$ arithmetic  operations \cite{BRITANAK200773}.
	
 Now we  conclude the current section with the following lemma on  the diagonalization of the matrix $B$ defined  in \eqref{importantmatrix}.
	\begin{lemma}\label{lemma2-1}
		It holds that
		\begin{align}\label{BSC_relation}
			S^{-1} BC= S B C^{-\tr}  = \Gamma,
		\end{align}
		where 
		\begin{align}\label{matrix_gamma}
			\Gamma {\color{black}:=} \left[ \bm 0, \Lambda\right]\in\RR^{(n-1)\times n},
		\end{align}
		and $\Lambda$ is defined as in \eqref{Lambda}.
	\end{lemma}
	
	\begin{proof}
		Firstly, one readily gets that $B\bm c_0 = \bm 0$.
		Then for $1\le j\le n-1,$ we have
		\begin{align*}
			B \bm c_j&= \left( \cos \tfrac{3j\pi}{2n} - \cos\tfrac{j\pi}{2n}, \cdots, \cos\tfrac{(2n-1)j\pi}{2n} - \cos\tfrac{(2n-3)j\pi}{2n}\right)^{\tr}
			\\
			&=-2\sin \tfrac{j \pi}{2n} \left( \sin \tfrac{j \pi}{n},  \sin \tfrac{2j \pi}{n}, \dots,
			\sin \tfrac{(n-1)j \pi}{n}\right)^{\tr}
			= -2\sin \tfrac{j \pi}{2n} \bm s_j.
		\end{align*}
		Thus, it is obtained that
		\begin{align}\label{BSC_equ1}
			B C & = \left(\bm s_1,\bm s_2,\cdots \bm s_{n-1} \right) [\bm 0,\Lambda]
			=S\Gamma.
		\end{align}
  {\color{black} Besides, according to \eqref{S_inv} and \eqref{C_inv}, we also obtain 
  \begin{align*}
  BC^{-\tr}=S^{-1}\Gamma,
  \end{align*}
  }which proves \eqref{BSC_relation}. 
	\end{proof}

	\section{Exact eigen-decompositions of the discrete Maxwell's eigenvalue problem}\label{exactsolution} 
	Before coming to the fast algorithm for solving \eqref{Kron-eq-01-sou}-\eqref{Kron-eq-02-sou},
	we explore the corresponding eigen-decomposition,
	\begin{subequations}
		\begin{align}\label{Kron-eq-01-x}
			UA-B^\tr VB^\tr =\dfrac{h^2}{6} \lambda^h  UD,
			\\
			\label{Kron-eq-01-y}
			-BUB+AV = \dfrac{h^2}{6} \lambda^h DV,
	\end{align}
	\end{subequations}
	which serves as a discrete version of the Maxwell's eigenvalue problem \eqref{Maxwell_eigen} with the boundary condition \eqref{essential_BCs}.

%
	To this end, we denote 
	\begin{align}\label{hat_UV}
		\widehat{U} =  C^\tr U S\in\RR^{n\times(n-1)},\quad
		\widehat{V} =  S V  C\in\RR^{(n-1)\times n}.
	\end{align}
	Then, premultiplying by $ C^{\tr}$ and postmultiplying  by $ S$ in  \eqref{Kron-eq-01-x},
	premultiplying  by $ S^{\tr}$  and postmultiplying  by $ C$ in \eqref{Kron-eq-01-y},
	we derive 	by using \eqref{S_inv}, \eqref{diag_AD}, \eqref{C_inv} and \eqref{BSC_relation}  that
	\begin{subequations}
	\begin{align}\label{Kron-eq-05-x}
			\widehat{U} \Lambda^2 
			- \Gamma^\tr \widehat{V} \Gamma^\tr =\dfrac{h^2}{6} \lambda^h \widehat{U} \Sigma,
			\\
			\label{Kron-eq-05-y}
			- \Gamma \widehat{U}\Gamma  +\Lambda^2 \widehat{V}=\dfrac{h^2}{6} \lambda^h \Sigma \widehat{V}.
	\end{align}
	\end{subequations}	
	Further, premultiplying \eqref{Kron-eq-05-x} by $\Gamma$, postmultiplying \eqref{Kron-eq-05-y} by $\Gamma^{\tr}$, {\color{black}using the fact $\Gamma\Gamma^\tr=\Lambda^2$} and  summing up the resultants, we find that
	\begin{align*}
		\lambda^h\big( \Gamma \widehat{U}\Sigma  +  \Sigma  \widehat{V}\Gamma^{\tr}\big) = \bm 0.
	\end{align*}
	
	To proceed, 
	let us first assume that $\lambda^h\neq 0,$ it  then yields  that
	\begin{align}\label{divfree_cond}
		\Gamma\widehat{U}\Sigma  +  \Sigma \widehat{V}\Gamma^{\tr}=\bm 0,
	\end{align}
	which is exactly the discrete divergence-free constraint \eqref{Kron-eq-02-sou} on the eigen-functions.
	In fact, by \eqref{diag_AD}, \eqref{BSC_relation} and \eqref{hat_UV}, we derive
	\begin{align*}
		\Gamma\widehat{U}\Sigma  +  \Sigma \widehat{V}\Gamma^{\tr} =&\,
		(SBC^{-\tr})( C^{\tr}U S)(S^{-1} D S) + (S D S^{-1})(S V C)( C^{-1} B^{\tr} S)
		\\
		=&\,S(BUD+ DVB^\tr)S,
	\end{align*}
	which states the equivalence of \eqref{divfree_cond} and the divergence-free constraint \eqref{Kron-eq-02-sou}.
	Therefore, all eigenfunctions associated with non-zero eigenvalues of \eqref{Kron-eq-01-x}- \eqref{Kron-eq-01-y} are exactly the discrete  divergence-free solutions.  
	In the sequel, we combine   \eqref{Kron-eq-05-x}-\eqref{Kron-eq-05-y} with \eqref{divfree_cond}, then find that 
	\begin{align}
		&\begin{cases}
			\widehat{U} \Lambda^2  \Sigma^{-1}+
			\Gamma^\tr \Sigma^{-1} \Gamma\widehat{U} =\dfrac{h^2}{6} \lambda^h \widehat{U}
			,\\[0.5em]
			\widehat{V}\Gamma^{\tr} \Sigma^{-1} \Gamma +\Sigma^{-1}\Lambda^2 \widehat{V}=\dfrac{h^2}{6} \lambda^h  \widehat{V}.
		\end{cases}\label{nonzero_system-3}
	\end{align}
	As a consequence, the eigen-solutions of \eqref{Kron-eq-05-x}-\eqref{Kron-eq-05-y} with  non-zero eigenvalues read,
	\begin{align*}
		&\Big(\lambda^h_{i,j},
		\widehat{U}_{i,j}, \; \widehat{V}_{i,j}\Big)  = \Big(
		\frac{6}{h^2} \Big(   \frac{\tau_i^2}{\sigma_i} + \frac{\tau_j^2}{\sigma_j}\Big),
		{\sigma_i \tau_j}  \bm e_{i+1}^{n} \otimes \bm e_{j}^{n-1}{}^\tr,\;  -{\sigma_j \tau_i}  \bm e_i^{n-1}\otimes\bm e_{j+1}^{n}{}^\tr\Big),
		\\
		&\qquad 0\le i,j \le n-1,\  (i,j)\neq (0,0). 
	\end{align*}

	Meanwhile, if $\lambda^h=0$,  \eqref{Kron-eq-05-x}-\eqref{Kron-eq-05-y} reduces to 
	\begin{align}
		\widehat{U} \Lambda^2  =
		\Gamma^\tr \widehat{V} \Gamma^\tr, \qquad 
		\Gamma \widehat{U}\Gamma  =\Lambda^2 \widehat{V}.
	\end{align}
	Then one readily concludes the corresponding curl-free eigen-solutions as
	\begin{align*}
		\Big(\lambda^h_0, \widehat {U}_{i,j}^0, \; \widehat {V}_{i,j}^0\Big) = 
		\Big(0,\; {\tau_i} \bm e_{i+1}^{n} \otimes \bm e_j^{n-1}{}^\tr,\;  {\tau_j}  \bm e_i^{n-1} \otimes \bm e_{j+1}^{n}{}^\tr\Big),\quad  1\le i,j \le n-1.
	\end{align*}
	Finally, by \eqref{hat_UV}, we summarize the eigen-solutions in the following concise form.
	
	\begin{theorem}\label{hat_solution}
		 {\color{black}The $2n(n-1)$ eigen-solutions to  \eqref{Kron-eq-01-x}- \eqref{Kron-eq-01-y} are of the following form:}
		\begin{itemize}
			\item  $n^2-1$ eigen-solutions with discrete divergence-free constraint  are given by
			\begin{align}
				\label{DFeigensolutions}
				\left( \lambda^h_{i,j}, U_{i,j},V_{i,j} \right) =&
				\left( \dfrac{6}{h^2}\Big( \dfrac{\tau_i^2}{\sigma_i} + \dfrac{\tau_j^2}{\sigma_j}\Big), \tau_j \sigma_i\bm c_i  \otimes\bm s_j^\tr, -\tau_i \sigma_j  \bm s_i  \otimes\bm c_j^\tr\right),  
				\\
				&\qquad 0\le i,j\le n-1,\ (i,j)\neq (0,0).\notag
			\end{align}

			\item $(n-1)^2$ eigen-solutions   with discrete curl-free constraint are given by
			\begin{align*}
					\left(\lambda^h_0, U^0_{i,j},V^0_{i,j}\right) = \left(0, \tau_i\bm c_i \otimes\bm s_j^\tr,\tau_j \bm s_i \otimes\bm c_j^\tr\right),
				\quad 1\le i,j \le n-1.
			\end{align*}
		\end{itemize}
	\end{theorem}
	
 It is worth noting that the eigen-decomposition  above actually serves as the  Helmholtz-Hodge decomposition of the  vector field space $U^0_h$ at the discrete level,
 \begin{align*} 
 U^0_h=  & {\color{black}span} \big\{   \tau_j\sigma_i  \bm s_{j}^{\tr}  N^{(1)} \bm c_i -  \tau_i\sigma_j  \bm c_{j}^{\tr}  N^{(2)} \bm s_i  :\  0\le i,j\le n-1,\, (i,j)\neq (0,0)   \big\}
 \\
 & \oplus\, {\color{black}span} \big\{   \tau_i  \bm s_{j}^{\tr}  N^{(1)} \bm c_i +  \tau_j  \bm c_{j}^{\tr}  N^{(2)} \bm s_i
		:\  1\le i,j\le n-1   \big\},
			\end{align*}
where    $N^{(1)} =\big ( \bm N^{(1)}_{i,j}(x,y) \big)_{1\le i\le n-1\atop 1\le j\le n} ,\,   N^{(2)} =\big ( \bm N^{(2)}_{i,j}(x,y) \big)_{1\le i\le n\atop 1\le j\le n-1}  $.

	%

	\section{Fast solver for  Maxwell's  source problems}\label{fastsolver}
	In this section, we design a fast solver for the Maxwell's  source problem \eqref{Maxwell_source}-\eqref{essential_BCs} based on  its discrete linear  equations \eqref{Kron-eq-01-sou}-\eqref{Kron-eq-02-sou}.  Intuitively,  the  solution $(U,V)$ of the discrete Maxwell's equation \eqref{Kron-eq-01-sou}-\eqref{Kron-eq-02-sou}
		can be expanded as a series of the divergence-free eigenfunctions $(U_{i,j},V_{i,j} )$, which can be  obtained through a derivation process similar to  that of \eqref{DFeigensolutions}.
	It then suffices to explore a fast algorithm for solving  $(\widehat U, \widehat V) = (C^{\tr}US, SVC)$ from
	the  following system derived from \eqref{Kron-eq-01-sou}-\eqref{Kron-eq-02-sou}:
	\begin{subequations}
	\begin{align}\label{Kron-eq-03-sou-x}
			&\widehat{U} \Big( \Lambda^2 +\dfrac{h^2}{6}\alpha \Sigma \Big) - \Gamma^\tr \widehat{V} \Gamma^\tr = h^2 \widehat{F},
			\\
			\label{Kron-eq-03-sou-y}
			&-\Gamma \widehat{U} \Gamma + \Big( \Lambda^2+\dfrac{h^2}{6}\alpha \Sigma\Big)\widehat{V} = h^2 \widehat{G},
		\\
		&\Gamma \widehat{U} \Sigma + \Sigma \widehat{V} \Gamma^\tr = \bm 0,
		\label{Kron-eq-03-sou-div}
	\end{align}
	\end{subequations}
	where  
	$
		\widehat{F} =  C^\tr F  S\in\RR^{n\times(n-1)}$ and $
		\widehat{G} =  S G  C\in\RR^{(n-1)\times n}.$

	Combining \eqref{Kron-eq-03-sou-x} and \eqref{Kron-eq-03-sou-div}, we obtain
	\begin{align}\label{solve-03}
		\widehat{U}
		\Big( \Lambda^2 +\frac{h^2}{6}\alpha\Sigma\Big)
		+\Gamma^{\tr} \Sigma^{-1}\Gamma \widehat{U}\Sigma = h^2 \widehat{F}.
	\end{align}
	Similarly,  we combine  \eqref{Kron-eq-03-sou-y} and \eqref{Kron-eq-03-sou-div}, and  get that
	\begin{align}\label{solve-04}
		\Sigma \widehat{V}  \Gamma^{\tr} \Sigma^{-1} \Gamma + 
		\Big( \Lambda^2 +\dfrac{h^2}{6}\alpha\Sigma\Big)
		\widehat{V} = h^2 \widehat{G}.
	\end{align}
{\color{black} The entries of \eqref{solve-03}-\eqref{solve-04}  are
\begin{align}
    \widehat{U}(i+1,j)\left(\tau_j^2+\frac{h^2}{6}\alpha\sigma_j+\frac{\tau_i^2}{\sigma_i}\sigma_j\right)=h^2\widehat{F}(i+1,j),\qquad 0\le i\le n-1, 1\le j\le n-1,\\
    \widehat{V}(i,j+1)\left(\frac{\tau_j^2}{\sigma_j}\sigma_i+\tau_i^2+\frac{h^2}{6}\alpha\sigma_i\right)=h^2\widehat{G}(i,j+1),\qquad 1\le i\le n-1, 0\le j\le n-1.
\end{align}

Hence,} we summarize the above computational  procedures into  Algorithm \ref{alg:divfree} with the complexity of $\mathcal{O}(n^2\log n)$.
	\begin{algorithm}[htb] 
		\caption{ Solve $U$ and $V$ in \eqref{Kron-eq-01-sou}-\eqref{Kron-eq-02-sou}.} 
		\label{alg:divfree} 
		\begin{algorithmic}[1] 
			\Require 
			the arrays $F$, $G$; the partition size $n$; the mesh size $h$; the coefficient $\alpha$.
			\Ensure 
			the solutions $U$ and $V$.
			\State Compute $\sigma_i$ and $\tau_i$ for $0\le i\le n-1$ by
			\begin{align*}
				\tau_i = -2\sin \frac{i\pi}{2n},\qquad 
				\sigma_i = 2\Big( 2+\cos\frac{i\pi}{n}\Big).
			\end{align*}
			\State Perform the FSTs/FCTs: 
						\begin{align*}
				\widehat{F}=C^{\tr}FS, \qquad \widehat{G}=SGC.
			\end{align*}
			\State Evaluate $\widehat{U}= \big[ \widehat{U}(i+1,j) \big]_{0\le i\le n-1\atop 1\le j\le n-1}$ and $\widehat{V}=\big[ \widehat{V}(i,j+1) \big]_{1\le i\le n-1\atop 0\le j\le n-1}$ by
			\begin{subequations}
			\begin{align}\label{essential_UV_hat_U}
					&\widehat{U}(i+1,j)= \dfrac{6h^2 \sigma_i\widehat{F}(i+1,j)}{6\tau_j^2\sigma_i+h^2\alpha\sigma_i\sigma_j +6\tau_i^2\sigma_j}
					=  \dfrac{h^2 \widehat{F}(i+1,j)}{(\lambda_{i,j}^h+\alpha) \frac{h^2}{6}\sigma_j},\\
					\label{essential_UV_hat_V}
					&\widehat{V}(i,j+1)=   \dfrac{6h^2 \sigma_j\widehat{G}(i,j+1)}{6\tau_i^2\sigma_j+h^2\alpha\sigma_i\sigma_j +6\tau_j^2\sigma_i}
					= \dfrac{h^2 \widehat{G}(i,j+1)}{(\lambda_{i,j}^h+\alpha) \frac{h^2}{6}\sigma_i}.
			\end{align}
			\end{subequations}
						\State Compute $U$ and $V$ by the inverse FSTs/FCTs:
			\begin{align*}
				U=
				\frac{4}{n^2}C\widehat{U}S, \qquad V=\frac{4}{n^2}S\widehat{V}C^{\tr}.
			\end{align*}
			\Return $U$ and $V$.
		\end{algorithmic} 
	\end{algorithm}

	\begin{remark}
		Premultiplying \eqref{Kron-eq-03-sou-x} by $\Gamma$  and postmultiplying \eqref{Kron-eq-03-sou-y}  by $\Gamma^\tr$ then summing up the results  yield
		\begin{align}\label{compatibility-matrix}
			\Gamma\widehat{F}+\widehat{G} \Gamma^\tr=\bm 0,
		\end{align}
		which are the equivalent compatibility conditions held by the right-hand side terms,
		\begin{align*}
			BF + G B^{\tr} = \bm 0.
		\end{align*}
		As a consequence, all columns except the first one of  $\widehat G$ are implicitly determined by $\widehat F$,   while  $\widehat G(:,1)$  still need to be obtained  through a direct calculation,
		\begin{align*}
			\widehat{G}(:,1) = SGC(:,1) = \tfrac{1}{\sqrt{2}} S \Big(\sum\limits_{k=1}^n G(j,k) \Big)_{j=1}^{n-1},
		\end{align*}
		with the complexity of  $3(n-1)^2$. Thus, $\widehat{V}$ shall be equivalently computed by
		\begin{align*}
			\widehat{V}(i,1) = \dfrac{6h^2\widehat{G}(i,1)}{6\tau_i^2+h^2\alpha\sigma_i}, 
			\widehat{V}(i,j+1) =  -\dfrac{\tau_i}{\tau_j}
			\dfrac{6h^2\sigma_j \widehat{F}(i+1,j)}{6\tau_i^2\sigma_j + h^2\alpha\sigma_i\sigma_j +6\tau_j^2\sigma_i},  1\le i,j \le n-1.
		\end{align*}
		
	\end{remark}

	\begin{remark}
	\label{Remark2}
		 It is worthy to note that  the Maxwell's equations  \eqref{Maxwell_source} and  \eqref{essential_BCs}  without the divergence equation \eqref{divergence}  are  still  well-posed for a proper parameter  $\alpha\neq 0$. 
		In fact, by eliminating the variable $\widehat{V}$ in  \eqref{Kron-eq-03-sou-x}-\eqref{Kron-eq-03-sou-y}, we  have
		\begin{align*}
			\widehat{U} \Big(\Lambda^2+\frac{h^2}{6} \alpha\Sigma\Big) 
			-\Gamma^{\tr}\Big(\Lambda^2+\frac{h^2}{6} \alpha\Sigma\Big)^{-1} \Gamma  \widehat{U} \Lambda^2
			= h^2   \Gamma^{\tr}\Big(\Lambda^2+\frac{h^2}{6} \alpha\Sigma\Big)^{-1}  \widehat{G} \Gamma^{\tr}
			+ h^2 \widehat{F}.
		\end{align*}
		Similarly, we find that
		\begin{align*}
			\Big( \Lambda^2 +\frac{h^2}{6}\alpha \Sigma\Big)\widehat{V} -\Lambda^2\widehat{V} \Gamma^\tr \Big( \Lambda^2 +\frac{h^2}{6}\alpha \Sigma\Big)^{-1}\Gamma = h^2\Gamma\widehat{F} \Big( \Lambda^2 +\frac{h^2}{6}\alpha \Sigma\Big)^{-1}\Gamma +h^2\widehat{G}.
		\end{align*}
		As a result, we conclude that
		\begin{subequations}
		\begin{align}
		\label{gauss_essential_UV_hat_U}
			\widehat{U}(i+1,j) = \frac{6h^2 \left( 6\tau_i^2 +h^2\alpha \sigma_i\right) \widehat{F}(i+1,j) + 36h^2\tau_i\tau_j\widehat{G}(i,j+1)}{ \left( 6\tau_i^2 +h^2\alpha \sigma_i\right)  \left( 6\tau_j^2 +h^2\alpha \sigma_j\right) - 36\tau_i^2\tau_j^2},\; \\0\le i\le n-1, \,1\le j\le n-1,\notag\\
			\label{gauss_essential_UV_hat_V}
			\widehat{V}(i,j+1) = \frac{6h^2 \left( 6\tau_j^2 +h^2\alpha \sigma_j\right) \widehat{G}(i,j+1) + 36h^2\tau_i\tau_j\widehat{F}(i+1,j)}{ \left( 6\tau_i^2 +h^2\alpha \sigma_i\right)  \left( 6\tau_j^2 +h^2\alpha \sigma_j\right) - 36\tau_i^2\tau_j^2},\; \\1\le i\le n-1, \,0\le j\le n-1 .\notag
		\end{align}
		\end{subequations}

		 Some emphases need to be  placed here.  For any general source term $\bm f$ in \eqref{Maxwell_source},  \eqref{gauss_essential_UV_hat_U}-\eqref{gauss_essential_UV_hat_V} hints $\frac{\alpha}{6}\big( \Gamma \widehat U \Sigma + \Sigma\widehat V \Gamma^\tr\big)  =  \Gamma \widehat F + \widehat  G \Gamma^\tr$ or its equivalence $\frac{\alpha}{6}\big(B U D + D V B^{\tr}\big) =B F+ G B^{\tr}  $, 
  indicating that our fast solver  inherently preserves	 
		 the Gauss's law $\nabla \cdot \bm u = \frac{1}{\alpha}\nabla \cdot \bm f:=\rho $  at the discrete level.
		
			In particular, if the source term $\bm f$ is divergence-free, then  the compatibility relation \eqref{compatibility-matrix} is satisfied, and			 the divergence-free constraint 
		 of the numerical solution computed from \eqref{gauss_essential_UV_hat_U}-\eqref{gauss_essential_UV_hat_V}  is guaranteed at the discrete level.
		This gives the equivalence between \eqref{gauss_essential_UV_hat_U}-\eqref{gauss_essential_UV_hat_V}  and  \eqref{essential_UV_hat_U}-\eqref{essential_UV_hat_V}.
	\end{remark}
	
	\begin{remark}
		Let us  revisit the fast direct solvers \cite{Zhang2015How, Averbuch1998A, Hockney1970The} designed for the scalar Helmholtz equation $-\Delta u+\alpha u = f$ in the rectangle  $[0,1]^2$ with the zero boundary conditions.  The discrete equation for the Helmholtz equation is derived using
	 the bilinear finite element method on a uniform mesh and can be written as follows:
		\begin{align}\label{matrix-poisson}
			AUD+DUA+ \frac{h^2}{6} \alpha DU D=6F, \quad A,D,U\in \RR^{(n-1)\times (n-1)}, \  h =\frac{1}{n},
		\end{align}
		where $U=\big[ u(i,j) \big]_{ 1\le i \le n-1\atop 1\le j\le n-1}$ with its entries being the nodal values of the finite element solution.
		
		Owing to  the eigendecompositions \eqref{diag_AD} together with \eqref{S_inv},  $U=S    \widehat U S$ can be calculated  through the transformed right-hand side $\widehat F:= S F S$ and 
		\begin{align*}
& \widehat U(i,j) =\frac{36 \widehat F(i,j)}{6\tau_i^2\sigma_j + 6\sigma_i\tau_j^2+h^2\alpha  \sigma_i \sigma_j}
 = \frac{ 6\widehat F(i,j)}{ ( \lambda_{i,j}^h+\alpha)  \frac{h^2}{6}\sigma_i \sigma_j},\quad 1\le i,j\le n-1,
 		\end{align*}
		where  we note that $\lambda_{i,j}^h=\frac{6}{h^2}\left(\frac{\tau_i^2}{\sigma_i}+\frac{\tau_j^2}{\sigma_j}\right)$ is also an eigenvalue of the corresponding  discrete Poisson eignevalue problem $AUD+DUA= \frac{h^2}{6} \lambda^h DU D$. 
		
		Noting that  the fast Helmholtz solvers based on \eqref{matrix-poisson} can be evaluated with an  arithmetic complexity of $\mathcal{O}(n^2 \log n)$
		 using FSTs/FCTs, we ultimately conclude that our fast Maxwell solver  fully rivals existing fast Helmholtz solvers in the literature. This capability enables the efficient solution of Maxwell's equations, akin to solving the Helmholtz equation.
	\end{remark}
	%
	%
	%


	\section{Fast solver for Maxwell's equations  with natural boundary conditions}\label{natural_bc}
	In this section, we shall generalize the fast solver to the source problem \eqref{Maxwell_source}-\eqref{divergence} with the  natural  boundary conditions
	\eqref{natural_BCs}.
	In  this scenario, we express the numerical solution $\bm u_h$ by global basis functions   as
	\begin{align}\label{uh_series_nat}
		\bm u_h = \sum\limits_{j=1}^{n} \sum\limits_{i=0}^{n}  \widehat{u}_{i,j}\bm N_{i,j}^{(1)}(x,y) + 
		\sum\limits_{j=0}^{n} \sum\limits_{i=1}^{n}   \widehat{v}_{i,j}\bm N_{i,j}^{(2)}(x,y)\in U_h,
	\end{align}
	and then
	obtain the numerical approximation scheme: to find $\bm u_h\in U_h $ such that
	\begin{subequations}
	\begin{align}\label{scheme_source_natural}
			&\left( {\color{black} \rot\ } \bm u_h , {\color{black} \rot\ } \bm w_h\right) +\alpha (\bm u_h,\bm w_h)  = \left(\bm f, \bm w_h\right),\quad &\forall \bm w_h\in U_h,\\
			\label{scheme_divergence_natural}
			&\left( \bm u_h,\nabla q_h\right) = 0,\quad &\forall q_h\in S_h.
	\end{align}
	\end{subequations}
	
	Combining \eqref{scheme_source_natural}, \eqref{scheme_divergence_natural}, and \eqref{uh_series_nat}, we obtain the linear algebraic system for source problem with natural boundary conditions and the divergence-free constraint:
	\begin{subequations}	
	\begin{align}
		&\begin{cases}
			U\bar{A} - \bar{B}V\bar{B} +\dfrac{\alpha h^2}{6} U \bar{D}  = h^2 F,\\
			-\bar{B}^\tr U \bar{B}^\tr +\bar{A} V + \dfrac{ \alpha h^2}{6} \bar{D}  V = h^2 G,
		\end{cases} \label{sou_natural_vec1} \\ 
		&\bar{B}^\tr U \bar{D}  +  \bar{D}  V\bar{B} = \bm 0,\label{sou_natural_vec2}
	\end{align}
		\end{subequations}
	where {\color{black}$U\! =\! \big( \widehat{u}_{i,j}\big)_{0\le i\le  n\atop 1\le j \le n}^\tr$, $V\! =\!\big( \widehat{v}_{i,j}\big)_{1\le i\le  n\atop 0\le j \le n}^\tr$,
	$F\!=\!\big( (\bm f,\bm N^{(1)}_{i,j}) \big)_{0\le i\le  n\atop 1\le j \le n}^\tr$,
	$G\!=\!\big( (\bm f,\bm N^{(2)}_{i,j}) \big)_{1\le i\le  n\atop 0\le j \le n}^\tr$.}
	The matrices $\bar{A}, \bar{D} \in\RR^{(n+1)\times (n+1)}$  and $\bar{B}\in\RR^{n\times (n+1)}$  are defined as
	{\addtolength{\arraycolsep}{-0.15em}
		\begin{align}\label{AD_tilde}
			\begin{split}
				&\bar A\!=\!
				\begin{pmatrix}
					1&-1 &\\
					-1 &2 &-1\\
					&\ddots&\ddots&\ddots\\
					&&-1&2&-1\\
					&&&-1&1
				\end{pmatrix},\  \bar  D\! =\!\begin{pmatrix}
					2&1 &\\
					1 &4 &1\\
					&\ddots&\ddots&\ddots\\
					&&1&4&1\\
					&&&1&2
				\end{pmatrix},\ 
				\bar B\!=\!\begin{pmatrix}
					\!-1\! &\!1\! &&&\\
					&\!-1\!&\!1\!&&\\
					& &\ddots&\ddots&\\
					&&&\!-1\!&\!1\!
				\end{pmatrix}.
			\end{split}
		\end{align}
	}

	We shall design the fast solver for \eqref{sou_natural_vec1}-\eqref{sou_natural_vec2}
	by exploring the exact discrete eigen-decompositions of the corresponding  Maxwell's eigenvalue problem  \eqref{Maxwell_eigen} with the natural boundary condition \eqref{natural_BCs}.
	
	As the computational procedures  resemble those in Section \ref{exactsolution} and Section \ref{fastsolver} for the case with essential boundary conditions, we  will present our conclusions in this section, omitting the detailed procedures.

	\subsection{Eigenvalue problem}
	\label{Eigen-natural}
	We first focus on the discretizations of \eqref{Maxwell_eigen} and  \eqref{natural_BCs}. 
	The corresponding discrete eigenvalue problem is: to find the non-trivial $(\bar \lambda^h,  U,   V)\in \RR\times \RR^{n\times(n+1)} \times \RR^{(n+1)\times n}$ such that
	\begin{align}\label{bar_eig_vec}
		\begin{cases}
		  U\bar{A} - \bar{B}  V \bar{B} = \dfrac{h^2}{6}\bar \lambda^h  U \bar D ,\\
			-\bar B^\tr  U \bar B^\tr +\bar A V = \dfrac{h^2}{6} \bar \lambda^h  \bar D  V.
		\end{cases}
	\end{align}
	
	To deal with problems associated with natural boundary conditions, 
	we need to introduce other kinds of DST and DCT.  
	For $n\in\NN,$ let 
	\begin{align}
		\label{DST_natural}
		& \bar S =  \left( \bar{\bm s}_1, \bar{\bm s}_2, \cdots, \bar{\bm s}_n\right),\qquad 
		\bar{\bm s}_j = \nu_j^n \Big( \sin \frac{j\pi}{2n}, \sin \frac{3 j\pi}{2n}, \dots, \sin \frac{(2n-1) j\pi}{2n}\Big)^\tr,
		\\
		\label{DCT_natural}
		&\bar C = \left( \bar{\bm c}_0,  \bar{\bm c}_1, \dots,  \bar{\bm c}_n\right) , 
		\qquad  \bar{\bm c}_j = \nu_j^n \Big( \frac{1}{\sqrt{2}}, \cos \frac{ j \pi}{n},   \cdots,  \cos \frac{(n-1)  j \pi}{n},\frac{(-1)^j}{\sqrt{2}}\Big)^\tr,
	\end{align}
	where $\nu_j^n$ is defined as in \eqref{nu_j}.
	Thus, multiplying a vector by $\bar{S}^\tr$, $\bar{S}$ and $\bar{C}$ represents a DST of the second type (DST-\IItwo), DST of the third type (DST-\IIthree) and DCT of the first type (DCT-\IIone) for the vector, respectively.
	Obviously, these two matrices satisfy \cite{BRITANAK200773}
	\begin{align}\label{SC_inv}
		\bar S^{-1} = \frac{2}{n} \bar S^\tr,\qquad \bar C^{-1} = \frac{2}{n} 
		\bar C .
	\end{align}
	
	%
	%
	
	Moreover, one readily checks that 
	\begin{align}\label{diag_ADtilde}
		\begin{split}
			&\bar C \bar M \bar A  \bar M \bar C^{-1} =   \bar C^{-1} \bar M \bar A  \bar M \bar C =   \frac{2}{n} \bar C \bar M \bar A  \bar M \bar C =   \bar{\Lambda}^2,
			\\
			& \bar C \bar M \bar D  \bar M \bar C^{-1} =   \bar C^{-1} \bar M \bar D  \bar M \bar C =\frac{2}{n}  \bar C \bar M \bar D  \bar M \bar C =   \bar{\Sigma},
		\end{split}
	\end{align}
	where
	\begin{align}
		&\bar M = \diag\big(\sqrt{2},1,\cdots,1,\sqrt{2}\big)\in\RR^{(n+1)\times(n+1)},\\
		&\label{Lambda_Sigma}
		\bar{\Lambda} = \diag(\tau_0,\tau_1, \cdots \tau_n),\quad
		\bar{\Sigma} = \diag( \sigma_0, \sigma_1,\cdots, \sigma_n),
	\end{align}
	and $\tau_i$ and $\sigma_i$ are defined as in \eqref{Lambda} and \eqref{Sigma} respectively. 
	Similar to the conclusions in Lemma \ref{lemma2-1}, we have the following lemma.

	\begin{lemma}
		It holds that 
		\begin{align}\label{B_diag}
			\bar  S^{-1} \bar B \bar M \bar C =  \bar  S^{\tr} \bar B \bar M \bar C^{-1}  =   \bar{\Gamma}, 
		\end{align}
		where
		\begin{align}
			\bar{\Gamma} = [\bm 0, \underline{\Lambda}]\in\RR^{n\times (n+1)},\quad 
			\underline{\Lambda} = \diag(\tau_1,\tau_2,\cdots, \tau_n).
		\end{align}
	\end{lemma}

	To solve \eqref{bar_eig_vec}, we first
	denote 
	\begin{align}\label{hat_UV1}
		\widehat{U} = \bar{S}^\tr  U \bar{M}^{-1} \bar{C}\in\RR^{n\times(n+1)},\quad \widehat{V} = \bar{C} \bar{M}^{-1}  V \bar{S} \in \RR^{(n+1)\times n}.
	\end{align}

	Applying a similar treatment to the computational procedures of \eqref{Kron-eq-01-x}- \eqref{Kron-eq-01-y}, we  obtain all eigensolutions of \eqref{bar_eig_vec},  providing an implementation of  the Hodge-Helmholtz decomposition for the approximation space $U_h$.
	\begin{theorem}\label{eig_solus_natural}
  {\color{black}The $2n(n+1)$ eigen-solutions to  \eqref{bar_eig_vec} are of the following form:}
		\begin{itemize}
			\item $n^2$ eigen-solutions with discrete divergence-free constraint are given by
			\begin{align*}
					(\bar \lambda_{i,j}^h ,  U_{i,j},  V_{i,j}  )= \bigg(\dfrac{6}{h^2} \Big( \dfrac{\tau_i^2}{\sigma_i} +\dfrac{\tau_j^2}{\sigma_j}\Big),  {\tau_j\sigma_i}  \bar{\bm s}_i  \otimes\tilde{\bm c}_j^\tr, -{\tau_i\sigma_j} \tilde{\bm c}_i  \otimes\bar{\bm s}_j^\tr\bigg);
\ 1\le i,j\le n.
			\end{align*}			
			
			\item $(n+1)^2-1$ eigen-solutions  { with discrete curl-free constraint} are given by
			\begin{align*}
					(\bar \lambda_0^h, U_{i,j}^0,  V_{i,j}^0) =  \left(0, {\tau_i}\bar{\bm s}_i  \otimes\tilde{\bm c}_j^\tr , {\tau_j} \tilde{\bm c}_i \otimes\bar{\bm s}_j^\tr\right),  \quad   0\le i,j\le n,\ (i,j)\neq (0,0),
			\end{align*}
		\end{itemize}
		where $ \tilde{\bm c}_j = M\bar {\bm c}_j= \nu_j^n\big(1, \cos\frac{j\pi}{n},   \cdots,  \cos \frac{(n-1)  j \pi}{n}, (-1)^j  \big)^{\tr}$ for $0\le j\le n.$
	\end{theorem}

	\subsection{Source problem}
	By introducing $(\widehat{U}, \widehat{V}) = (\bar{S}^\tr U \bar{M}^{-1} \bar{C}, \bar{C} \bar{M}^{-1} V \bar{S})$, the system
	\eqref{sou_natural_vec1}-\eqref{sou_natural_vec2}
	is equivalent to the following equations
	\begin{subequations}
	\begin{align}\label{natu_source_sys-x}
			&\widehat{U}\Big( \bar{\Lambda}^2 +\dfrac{h^2}{6}\alpha \bar{\Sigma}\Big) - \bar{\Gamma} \widehat{V} \bar{\Gamma}  = h^2\widehat{F},
			\\
			\label{natu_source_sys-y}
			&-\bar{\Gamma}^\tr \widehat{U} \bar{\Gamma}^\tr+ \Big( \bar{\Lambda}^2+\dfrac{h^2}{6}\alpha \bar{\Sigma}\Big) \widehat{V} = h^2\widehat{G},
			\\
			\label{natu_div_sys}
			&\bar{\Gamma}^\tr \widehat{U}\bar{\Sigma} + \bar{\Sigma} \widehat{V} \bar{\Gamma} = \bm 0,
	\end{align}
	\end{subequations}
	where 
	$
		\widehat{F}= \bar S^\tr F\bar{M} \bar C\in \RR^{n\times(n+1)}$ and $ \widehat{G} =  \bar C \bar{M}G  \bar S\in\RR^{(n+1)\times n}$.

	From the divergence equation \eqref{natu_div_sys}, we first conclude that
	\begin{align}\label{soluUV0}
		\widehat{U}(:,1) = \bm 0,\quad \widehat{V}(1,:) = \bm 0.
	\end{align}
	Then eliminating $\widehat{V}$ from  \eqref{natu_source_sys-x} by using \eqref{natu_div_sys} gives
	\begin{align}\label{solve-05}
		\widehat{U}\Big( \bar{\Lambda}^2 +\dfrac{h^2}{6}\alpha \bar{\Sigma}\Big) +\bar{\Gamma} \bar{\Sigma}^{-1} \bar{\Gamma}^\tr \widehat{U}\bar{\Sigma}  = h^2\widehat{F}.
	\end{align}
	Analogously, eliminating $\widehat{U}$ from  \eqref{natu_source_sys-y} by using \eqref{natu_div_sys} yields
	\begin{align}\label{solve-06}
		\bar{\Sigma}\widehat{V} \bar{\Gamma}\bar{\Sigma}^{-1} \bar{\Gamma}^\tr+ \Big( \bar{\Lambda}^2+\dfrac{h^2}{6}\alpha \bar{\Sigma}\Big) \widehat{V} = h^2\widehat{G}.
	\end{align}
	Thus,
	an efficient algorithm to solve \eqref{sou_natural_vec1}-\eqref{sou_natural_vec2} with the complexity of $\mathcal{O}(n^2\log n)$ is given in Algorithm \ref{alg:divfree_natural}.
	\begin{algorithm}[htb] 
		\caption{ Solve $U$ and $V$ in \eqref{sou_natural_vec1}-\eqref{sou_natural_vec2}.} 
		\label{alg:divfree_natural} 
		\begin{algorithmic}[1] 
			\Require 
			the arrays $F$, $G,$ $\bar{M}$; the partition size $n$; the mesh size $h$; the coefficient $\alpha$.
			\Ensure 
			the solutions $U$ and $V$.
			\State Compute $\sigma_i$ and $\tau_i$ for $0\le i\le n-1$ by
			\begin{align*}
				\tau_i = -2\sin \frac{i\pi}{2n},\qquad 
				\sigma_i = 2\Big( 2+\cos\frac{i\pi}{n}\Big).
			\end{align*}
			\State Perform the FSTs/FCTs:
			\begin{align*}
				\widehat{F}= \bar S^\tr F\bar{M} \bar C,\quad
				\widehat{G} =  \bar C \bar{M}G  \bar S.
			\end{align*}
			\State Evaluate $\widehat{U}=\big[ \widehat{U}(i,j)\big]_{1\le i\le n\atop 1\le j\le n+1} $ and $\widehat{V}=\big[ \widehat{V}(i,j)\big]_{1\le i\le n+1\atop 1\le j\le n}$ by
			\begin{subequations}
			\begin{align}\label{natural_UV_hat-U}
					&\widehat{U}(i,1) = 0,\
					\widehat{U}(i,j+1) = \frac{6h^2\sigma_i\widehat{F}(i,j+1) }{6\sigma_i\tau_j^2 +h^2\alpha\sigma_i\sigma_j +6\sigma_j \tau_i^2}
					=\dfrac{h^2 \widehat{F}(i,j+1)}{(\bar \lambda_{i,j}^h+\alpha) \frac{h^2}{6}\sigma_j},\\
					\label{natural_UV_hat-V}
					&\widehat{V}(1,j) = 0,\ \widehat{V}(i+1,j) = \frac{6h^2\sigma_j \widehat{G}(i+1,j)}{6\sigma_i\tau_j^2 +h^2\alpha\sigma_i\sigma_j +6\sigma_j \tau_i^2}=\dfrac{h^2 \widehat{G}(i+1,j)}{(\bar \lambda_{i,j}^h+\alpha) \frac{h^2}{6}\sigma_i}.
			\end{align}
			\end{subequations}
			\State Compute $U$ and $V$ by the inverse FSTs/FCTs:
			\begin{align}\label{UV_natural}
				U = \frac{4}{n^2} \bar S\widehat{U} \bar C \bar{M},\quad
				V = \frac{4}{n^2} \bar{M} \bar C \widehat{V}  \bar S^\tr.
			\end{align}
			\Return $U$ and $V$.
		\end{algorithmic} 
	\end{algorithm}
	\begin{remark}
		 Premultiplying  \eqref{natu_source_sys-x} by $\bar{\Gamma}^\tr$ and postmultiplying  \eqref{natu_source_sys-y}  by $\bar{\Gamma}$, we have
		\begin{align*}
			\bar{\Gamma}^\tr \widehat{F} + \widehat{G} \bar{\Gamma}=\bm 0,
		\end{align*}
		which is acually the compatibility conditions held by the right-hand side terms,
		\begin{align}
			\bar{B}^\tr F+ G\bar{B} = \bm 0.
		\end{align}
		Thus,  the array ${G}$ is actually not required to be  stored  in Algorithm \ref{alg:divfree_natural} and $\widehat{V}$ can be evaluated by
		\begin{align*}
			\widehat{V}(1,j) = 0,\quad\widehat{V}(i+1,j) = -\frac{\tau_i}{\tau_j} \frac{6h^2\sigma_j \widehat{F}(i,j+1)}{6\sigma_i\tau_j^2 +h^2\alpha\sigma_i\sigma_j +6\sigma_j \tau_i^2},\quad 1\le i,j\le n.
		\end{align*}
	\end{remark}

	\begin{remark}
		When $\alpha\ne 0$ is  a proper parameter, the Maxwell's equations
			\eqref{Maxwell_source} with the natural boundary conditions \eqref{natural_BCs} can also be  computed efficiently. Similar to 
			Remark \ref{Remark2}, if
			\eqref{natural_UV_hat-U}-\eqref{natural_UV_hat-V}  in Algorithm \ref{alg:divfree_natural}
			 are replaced by
		\begin{align*}
			\begin{split}
				&\widehat{U}(i,j+1)\!  =\! \dfrac{ 6h^2\left( 6\tau_i^2\! +\! h^2\alpha\sigma_i  \right)\widehat{F}(i,j+1) \! +\! 36h^2 \tau_i\tau_j \widehat{G}(i+1,j)  }{\left( 6\tau_i^2\! +\! h^2\alpha\sigma_i\right)\left( 6\tau_j^2\! +\! h^2\alpha\sigma_j\right) \! -\!  36\tau_i^2\tau_j^2 }, \; 1\! \le\!  i\! \le\!  n, \,0\! \le\!  j\! \le\!  n,\\
				&\widehat{V}(i+1,j) \! =\!  \dfrac{6h^2\left(6 \tau_j^2\! +\! h^2\alpha\sigma_j\right)\widehat{G}(i+1,j) \! +\! 36h^2\tau_i\tau_j \widehat{F}(i,j+1)  }{ \left(6 \tau_i^2\! +\! h^2\alpha\sigma_i\right)\left( 6\tau_j^2\! +\! h^2\alpha\sigma_j\right) \! - \! 36\tau_i^2\tau_j^2},\;\, 0\! \le\!  i\! \le \! n,\, 1\! \le\!  j\! \le\!  n,
			\end{split}
		\end{align*}
		then we design 
			a fast solver for the Maxwell's equations \eqref{Maxwell_source} and \eqref{natural_BCs}, which preserves the Gauss's law $\nabla \cdot \bm u = \frac{1}{\alpha}\nabla \cdot \bm f:=\rho $ intrinsically at the discrete level.
	\end{remark}
	%
	%
	
	\section{Fast solver for general Maxwell's equations}\label{solver_pressure}
	In this section, the fast solver is extended to the Maxwell's equation \eqref{Maxwell_source_div} under the Gauss's law \eqref{Maxwell_Gauss_law} with the boundary condition \eqref{Maxwell_source_boundary}.
	A mixed  finite element approximation scheme reads: to find $(\bm{u}_h, p_h)\in U_h^0\times S_h^0$ such that
		\begin{subequations}
	\begin{align}
	\label{mix_source_scheme_div}
			&\left({\color{black} \rot\ } \bm{u}_h, {\color{black} \rot\ } \bm{w}_h\right) +\alpha(\bm u_h,\bm w_h)+(\nabla p_h, \bm w_h)= (\bm{f},\bm{w}_h),\quad &\forall \bm{w}_h\in U_h^0,\\
			\label{mix_source_scheme_gauss}
			&(\bm{u}_h, \nabla q_h)=-(\rho,  q_h),\quad &\forall q_h\in S_h^0.
	\end{align}
	\end{subequations}
	
	Similar to the derivation of \eqref{Kron-eq-01-sou}-\eqref{Kron-eq-02-sou}, we obtain the linear algebraic system for \eqref{mix_source_scheme_div}-\eqref{mix_source_scheme_gauss}: to find  $(U,V,P)\in \RR^{n\times(n-1)}\times \RR^{(n-1)\times n }\times\RR^{(n-1)\times(n-1)}$ such that
	\begin{subequations}
	\begin{align}\label{Kron-eq-01-sou_div-x}
			&UA-B^\tr VB^\tr +\dfrac{h^2}{6} \alpha UD-\dfrac{h^2}{6} B^\tr P D= h^2 F,\\
			\label{Kron-eq-01-sou_div-y}
			&-BUB+AV +\dfrac{h^2}{6} \alpha DV-\dfrac{h^2}{6} DPB=  h^2 G,\\
			\label{Kron-eq-01-sou_div-g}
			&BUD+DVB^\tr=6H,
	\end{align}
	\end{subequations}
	where 
	$
		H=\left( (\rho, \varphi_{i,j})\right)_{1\le i,j\le  n-1}$ with 
		${\varphi}_{i,j}(x,y)$ being the global basis of $S^0_h$ at the nodal point $(ih,jh)$.

	Next, we shall explore how to solve \eqref{Kron-eq-01-sou_div-x}-\eqref{Kron-eq-01-sou_div-g} efficiently. Using  \eqref{hat_UV} and further introducing
	\begin{align}\label{hat_FG_div}
		\begin{split}
			 \widehat{F} =  C^\tr F  S,\quad
		\widehat{G} =  S G  C,\quad\widehat{P}=SPS,\quad  \widehat{H} =  SHS,
		\end{split}
	\end{align}
	one readily checks that  the following equivalent equations to \eqref{Kron-eq-01-sou_div-x}-\eqref{Kron-eq-01-sou_div-g},
	\begin{subequations}
	\begin{align}\label{hat_source_nat-U}
			&\widehat{U}\Big ( \Lambda^2+ \dfrac{h^2}{6}\alpha  \Sigma\Big)-\Gamma^\tr \widehat{V} \Gamma^\tr   -\dfrac{h^2}{6}\Gamma^\tr \widehat{P} \Sigma = h^2 \widehat{F},\\
			\label{hat_source_nat-V}
			&-\Gamma\widehat{U}\Gamma +\Big ( \Lambda^2 +\dfrac{h^2}{6}\alpha \Sigma\Big)\widehat{V}  -\dfrac{h^2}{6} \Sigma\widehat{P}\Gamma = h^2 \widehat{G},\\
			\label{hat_source_nat-g}
			&\Gamma\widehat{U} \Sigma +\Sigma \widehat{V} \Gamma^\tr = 6\widehat{H}.
	\end{align}
	\end{subequations}
	Firstly, we premultiply \eqref{hat_source_nat-U} by $\Gamma$, postmultiply \eqref{hat_source_nat-V} by $\Gamma^\tr$,
	 and then  find by summing up the resultants that 
	\begin{align}\label{solve-P}
		\Lambda^2\widehat{P}\Sigma+\Sigma\widehat{P}\Lambda^2=6( \alpha \widehat{H} - \Gamma  \widehat{F}-\widehat{G} \Gamma^\tr  ),
	\end{align}
	which gives  a diagonal solver for $\widehat P$.
	Next, combining  \eqref{hat_source_nat-U} and  \eqref{hat_source_nat-g} yields
		\begin{align*}
			\widehat{U}\Big ( \Lambda^2+ \dfrac{h^2}{6}\alpha  \Sigma\Big)-\Gamma^\tr (6\Sigma^{-1}\widehat{H}  - \Sigma^{-1}\Gamma\widehat{U} \Sigma)  -\dfrac{h^2}{6}\Gamma^\tr \widehat{P} \Sigma = h^2 \widehat{F},
		\end{align*}
		which is equivalent to 
		\begin{align}\label{solve-U}
			\widehat{U}\Big ( \Lambda^2+ \dfrac{h^2}{6}\alpha  \Sigma\Big) + \Gamma^\tr\Sigma^{-1}\Gamma\widehat{U} \Sigma =6\Gamma^\tr\Sigma^{-1}\widehat{H} +\dfrac{h^2}{6}\Gamma^\tr \widehat{P} \Sigma +  h^2 \widehat{F}.
		\end{align}
		Analogously, we have
		\begin{align}\label{solve-V}
			\Sigma\widehat{V}\Gamma^\tr\Sigma^{-1}\Gamma  +\Big ( \Lambda^2 +\dfrac{h^2}{6}\alpha \Sigma\Big)\widehat{V}  =6\widehat{H}\Sigma^{-1} \Gamma+\dfrac{h^2}{6} \Sigma\widehat{P}\Gamma +  h^2 \widehat{G}.
		\end{align}
	A fast solver for \eqref{Kron-eq-01-sou_div-x}-\eqref{Kron-eq-01-sou_div-g} is then presented
	in Algorithm \ref{alg:div}.	\begin{algorithm}[htb] 
		\caption{ Solve $U$, $V$ and $P$ in \eqref{Kron-eq-01-sou_div-x}-\eqref{Kron-eq-01-sou_div-g}.} 
		\label{alg:div} 
		\begin{algorithmic}[1] 
			\Require 
			the arrays $F$, $G$, $H$; the partition size $n$; the mesh size $h$; the coefficient $\alpha$.
			\Ensure 
			the solutions $U$, $V$ and $P$.
			\State Compute $\sigma_i$ and $\tau_i$ for $0\le i\le n-1$ by
			\begin{align*}
				\tau_i = -2\sin \frac{i\pi}{2n},\qquad 
				\sigma_i = 2\Big( 2+\cos\frac{i\pi}{n}\Big).
			\end{align*}
			\State Perform the FSTs/FCTs:  
						\begin{align*}
				\widehat{F}=C^{\tr}FS, \qquad \widehat{G}=SGC,\qquad \widehat{H} = SHS.
			\end{align*} 
			\State {\color{black}Evaluate $\widehat{P}=\big[\widehat{P}(i,j)\big]_{1\le i \le n-1\atop 1\le j\le n-1}$, $\widehat{U}=\big[\widehat{U}(i,j)\big]_{0\le i \le n-1\atop 1\le j\le n-1} $,  and
			$\widehat{V}=\big[\widehat{V}(i,j)\big]_{1\le i \le n-1\atop 0\le j\le n-1}$   by }
		\begin{subequations}
			\begin{align} \label{UV_solu_div-P}
					&{\color{black}\widehat{P}(i,j) = \dfrac{6\big(\alpha\widehat{H}(i,j) -\tau_i\widehat{F}(i+1,j)-\tau_j\widehat{G}(i,j+1)\big)}{\tau_i^2\sigma_j+\sigma_i\tau_j^2},}\\
                      \label{UV_solu_div-U}
						&\widehat{U}(i+1,j) 
						=\dfrac{36\tau_i\widehat H(i,j)+h^2\tau_i\sigma_i\sigma_j\widehat{P}(i,j)+6h^2\sigma_i\widehat{F}(i+1,j)}{6\sigma_i\tau_j^2+h^2\alpha\sigma_i\sigma_j+ 6\sigma_j\tau_i^2}, \\
						\label{UV_solu_div-V}
						&\widehat{V}(i,j+1) = 
						\dfrac{36\tau_j\widehat H(i,j)+h^2\tau_j\sigma_i\sigma_j\widehat{P}(i,j)+6h^2\sigma_j\widehat{G}(i,j+1)}{6\sigma_i\tau_j^2+h^2\alpha\sigma_i\sigma_j+ 6\sigma_j\tau_i^2}.
			\end{align}
			\end{subequations} 
			\State Compute $U$,  $V$ and $P$ by the inverse  FSTs/FCTs: 
			\begin{align*}
				U = \frac{4}{n^2}  C\widehat{U}  S,\quad
				V = \frac{4} {n^2} S  \widehat{V}  C^\tr,\quad P=\frac{4}{n^2}  S  \widehat{P}  S.
			\end{align*}
			\Return $U$, $V$ and $P$.
		\end{algorithmic} 
	\end{algorithm}
	\begin{remark} The fast solver described in Algorithm \ref{alg:div} preserves the Gauss's law at the discrete level.
		
		Further if the right-hand side term in \eqref{Maxwell_source_div} satisfies
		$
			\nabla \cdot \bm f =\alpha \nabla\cdot \bm u=\alpha\rho,
		$
		then the solution $p$ (resp. $P$) in \eqref{Maxwell_source_div} (resp. in Algorithm \ref{alg:div})  is identically   zero. 
		As a result,  the  compatibility condition  $\alpha \widehat{H}-	\Gamma  \widehat F-\widehat G\Gamma^{\tr}=0$ on right-hand side terms holds.
		This indicates that $$\frac{\alpha}{6}\big( \Gamma \widehat U \Sigma + \Sigma\widehat V \Gamma\big)  =  \Gamma \widehat F + \widehat  G \Gamma  \text{ or  equivalently } \frac{\alpha}{6}\big(B U D + D V B^{\tr}\big) =B F+ G B^{\tr}.  $$
		Then  one finds that Algorithm \ref{alg:div}  with  general $\rho$  (resp. $\rho=0$) is  degenerated into that in Remark 
		\ref{Remark2}  (resp.  Algorithm \ref{alg:divfree}) directly. 

		 Finally, we note that  a fast solver for  \eqref{Maxwell_source_div} with natural boundary conditions can be deduced analogously.
		 However, due to space limitations, we omit the details.
	\end{remark}

	\section{Numerical experiments}\label{example}
	To illustrate the {\color{black}effectiveness and efficiency} of the fast solver designed for the Maxwell's source equations, we conduct some numerical experiments in this section.
	We set $\Omega=(0,1)^2$ in all {\color{black} examples}.
	$n_x$ and $n_y$ are the partition size in the $x$- and $y$-direction respectively.
	When they are equal, we let $n=n_x=n_y$.
	Denote the error $\bm e_h = \bm{u}-\bm{u}_h.$ Theoretically, 
	the vector field $\bm u_h$ converges to $\bm u$ at the order of $\mathcal{O}(h)$ both  in $\bm L^2$-norm (i.e. $\|\bm e_h\|$) and {\color{black}$H(\tc;\Omega)$-semi norm} (i.e. $\|{\color{black} \rot\ } \bm{e}_h\|$).
	 First, we compare the performances of the fast Maxwell's solver (Fast Solver) 
	 with that of  {\color{black}LU factorization (LU Solver)}.
	To implement the fast solver, FFTW library is used to compute DSTs/DCTs. For {\color{black}LU solver}, 
	the resulting linear algebraic system is solved by the {\tt pardiso} from Intel MKL.
	The Intel C compiler ({\tt{icc}}) is used in both methods.
	Additionally, we demonstrate that the discrete matrix generated by the Maxwell's equation with constant coefficients  serves as an effective preconditioner for solving
	the equation with variable coefficients  using iterative methods. 
	We call the function {\tt dcg} in Intel MKL to implement the conjugate gradient/pre-conditioned conjugate gradient (CG/PCG) method for this problem. 
 
 {\color{black}The errors and convergence orders obtained from the fast solver for the initial terms are the same as those from the LU solver, so only the fast solver’s results are recorded for brevity.

}

	\begin{example}\label{exam1}
		Consider the Maxwell's equation \eqref{Maxwell_source}
		when $\alpha=0$ with {\color{black}non-homogenous} essential boundary conditions.
		The right-hand side function $\bm{f}$ is chosen such that the exact solution is
		$$\bm{u} =\big(3 x^3y^2+1,-3x^2y^3+3\big)^{\tr}.$$		
	\end{example}
	\begin{table}[htb]
		\caption{\small{Errors, convergence rates and computational time of {\color{black}LU solver} and fast solver in Example \ref{exam1}.
  }
  }\label{tab_exam1-1}
		\tabcolsep =0.3em
		\centerline{\small{
				\begin{tabular}{c|ccccc||cccc|c}
					\hline
					\quad & \multicolumn{4}{c}{ Errors and convergence orders} & &\multicolumn{5}{c} {Computational time}\\
					\hline
					\quad & \multicolumn{4}{c}{ Fast Solver} & &\multicolumn{3}{c} {{\color{black}LU Solver}}&& Fast Solver\\
					\hline
					$n$ & $\|\bm e_h \|$ & order & $\|{\color{black} \rot\ }\bm e_h \|$ & order& & {\color{black}factor.}(s) & {\color{black}solution}(s)& total(s) && total(s)\\
					\hline
					$2^7$ &5.77e-03& &2.05e-02& && 0.286 & 0.024& 0.310 && 0.00068\\
					$2^8$ &2.88e-03& 1.00 &1.02e-02& 1.00 &&1.344 & 0.128& 1.472 && 0.00306\\
					$2^9$ &1.44e-03& 1.00 &5.11e-03& 1.00  &&6.684 & 0.558& 7.242&& 0.01466\\
					$2^{10}$& 7.18e-04	&1.00 & 2.56e-03& 1.00  &&  32.958&  2.408 &	35.366&& 0.06000 \\
					$2^{11}$ & 3.59e-04& 1.00 & 1.28e-03& 1.00  && 169.610 & 10.325 &  179.935 && 0.33000\\
					$2^{12}$& 1.79e-04 & 1.00 & 6.39e-04& 1.00  && 953.470 & 44.230 & 997.700 && 1.63000 \\
					$2^{13}$& 8.97e-05 & 1.00 & 3.20e-04	& 1.00  && -	&-	&-	&&8.54000 \\
					$2^{14}$& 4.49e-05	& 1.00 & 1.60e-04& 1.00  && -	&-	&-	&&41.3400 \\
					$2^{15}$& 2.24e-05& 1.00 & 7.99e-05& 1.00  && -	&-	&-	&&178.550 \\
					\hline
		\end{tabular}}}
	\end{table}

	 Given the divergence-free constraint ensured
	 by a Lagrange multiplier when $\alpha=0$, we derive a saddle-point problem from the discretization \eqref{mix_source_scheme_div}-\eqref{mix_source_scheme_gauss} with $\rho=0$.
	Table \ref{tab_exam1-1} presents the corresponding errors, convergence orders, and computational  times for both the {\color{black}LU solver} and the  fast solver. 
	 The corresponding convergence orders   align with the expected theoretical orders. 
	The CPU time (in seconds) is also provided.
	For  the {\color{black}LU solver}, we record factorization phase time, solution phase time, and total time for a direct  solving process. It is evident that the factorization phase consumes the majority of the time.
 In contrast, the fast solver  significantly improves computational time by orders of magnitude.

\begin{figure}[htb]
		\begin{minipage}[t]{0.45\linewidth}
			\hspace*{\fill}\includegraphics[width=\textwidth]{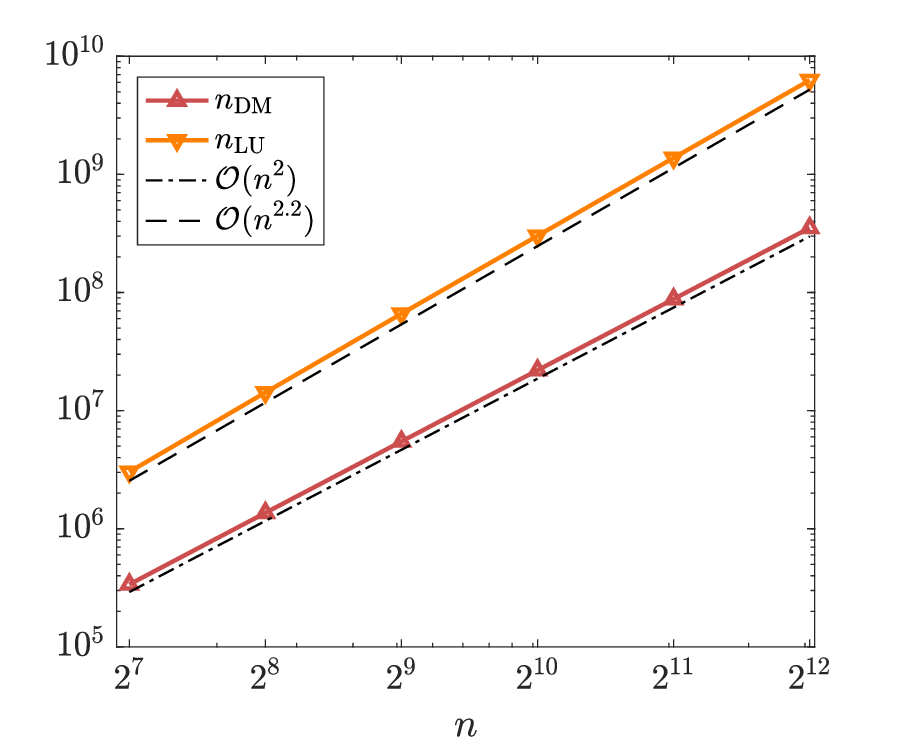}
			\vspace{-2em}
			\caption{\small{The number of non-zeros $n_{\rm DM}$ and $n_{\rm LU}$ versus $n$  in Example \ref{exam1}.}}\label{fig_exam1-1}
		\end{minipage}\hfill%
		\begin{minipage}[t]{0.45\linewidth}
			\centering
			\includegraphics[width=\textwidth]{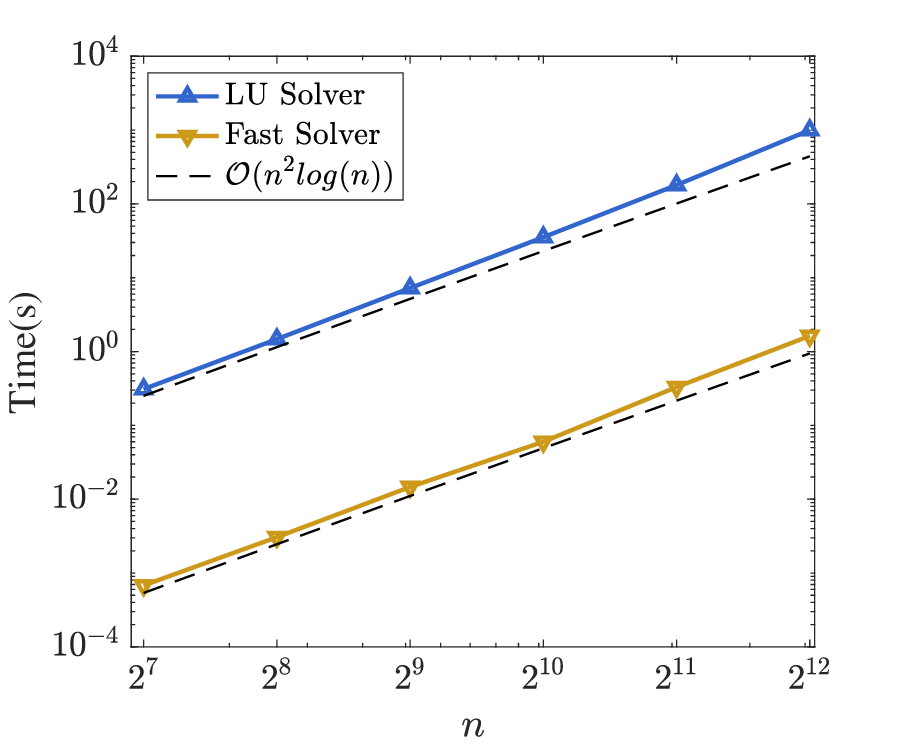}
			\vspace{-2em}
			\caption{\small{Total computational time versus $n$  in Example \ref{exam1}.}}\label{fig_exam1-2}
		\end{minipage}\hspace*{\fill}
	\end{figure}

	However, as $n$ becomes larger,  the {\color{black}LU solver} fails to derive the correct solution due to the  limitations of the available computational  resources (a workstation with 256GB memory). 
	Although the discrete matrix is highly sparse, 
	there is non-zero fill-in during matrix factorizations, requiring extra storage  to implement the {\color{black}LU solver}.
	From Figure \ref{fig_exam1-1}, we conclude that the numbers of non-zeros in the discrete matrix $(n_{\rm DM})$ and in the LU factorization $(n_{\rm LU})$ grow in the order of $\mathcal{O}(n^2)$  and $\mathcal{O}(n^{2.2})$, respectively. 
	It also suggests that the $n_{\rm LU}$  is much larger than $n_{\rm DM}$.
	We also  record the memory for the {\color{black}LU solver} and fast solver in Table \ref{mem_exam1-1}. When $n=6000,$  even considering the matrix factorizations,  the {\color{black}LU solver} requires 78.2\%  of the total memory, while  the fast solver only needs 0.2\%.

	\begin{table}[htb]
		\caption{\small{Memory of {\color{black}LU solver} and fast solver in Example \ref{exam1}. (total memory: 256GB) }}\label{mem_exam1-1}
		\centerline{\small{
				\begin{tabular}{c|ccccccccc}
					\hline
					$n$ & $2^7$ & $2^8$  &$2^9$  & $2^{10}$ & $2^{11}$ &  $2^{12}$ &6000 \\
					\hline
					{\color{black}mem(Fast Solver)} (MB)&1.5& 3.3&9.7&34 &131& 515 & 554\\
					{\color{black}mem(LU Solver)} (MB) &-& 262 &1049& 4719&20709&88343& 204997 \\
     {\color{black}mem(Fast Solver)/mem(LU Solver)}  (\%)&-&1.26&0.92&0.72&0.63&0.58&0.27\\
     \hline
		\end{tabular}}}
	\end{table}

	The asymptotic behaviors of computational time for these two different approaches are illustrated in Figure \ref{fig_exam1-2}. The time of  the fast solver grows in the order of $\mathcal{O}(n\!^2\!\log \!n)$.  In contrast,  the direct solver takes more time to solve, and the order grows faster than $\mathcal{O}(n^2\log n)$.
	
	\begin{example}\label{exam2}
		Consider the {\color{black} time-harmonic} Maxwell's equation \eqref{Maxwell_source} 
		when {\color{black}$\alpha=-1$} with homogeneous essential boundary conditions \eqref{essential_BCs}. The right-hand side function $\bm{f}$ is chosen such that the exact solution is
		$$\bm{u} = \big(\cos(\pi x)\sin(\pi y), -\sin(\pi x) \cos(\pi y)\big)^{\tr}.$$		
	\end{example}

 \begin{table}[htb]
{\color{black}		\caption{\small{Errors, convergence rates and computational time of {\color{black}LU solver} and fast solver in Example \ref{exam2}. }}\label{tab_exam2-1}
		\vspace{0.1cm}		
		\tabcolsep =0.3em
		\centerline{\small{
				\begin{tabular}{c|ccccc||cccc|c}
					\hline
					\quad & \multicolumn{4}{c}{ Errors and convergence orders} & &\multicolumn{5}{c} {Computational time}\\
					\hline
					\quad & \multicolumn{4}{c}{  Fast Solver} & &\multicolumn{3}{c} {{\color{black}LU Solver}}&& Fast Solver\\
					\hline
					$(n_x, n_y)$ & $\|\bm e_h \|$ & order & $\|{\color{black} \rot\ }\bm e_h \|$ & order& & {\color{black}factor.}(s) & {\color{black}solution}(s)& total(s) && total(s)\\
					\hline
					$(2^6, 2^7)$ &7.92e-03&  &4.98e-02& && 0.050 & 0.000 & 0.050  && 0.00031\\
					$(2^7, 2^8)$ &3.96e-03& 1.00 &2.49e-02& 1.00 &&0.240 & 0.010 & 0.250  &&0.00132 \\
					$(2^8, 2^9)$ &1.98e-03& 1.00 &1.24e-02& 1.00  &&1.190 	& 0.140 & 1.330 && 0.00669 \\
					$(2^9, 2^{10})$& 9.90e-04& 1.00 & 6.22e-03	&1.00 &&  5.660 & 0.590 & 6.250 && 0.02832 \\
					$(2^{10}, 2^{11})$ & 4.95e-04& 1.00 & 3.11e-03& 1.00  &&26.360& 2.480 & 28.840  && 0.13293\\
					$(2^{11}, 2^{12})$& 2.48e-04& 1.00 & 1.56e-03& 1.00  && 118.370& 8.880 & 127.250 && 0.62570  \\
					$(2^{12}, 2^{13})$& 1.24e-04& 1.00 & 7.78e-04& 1.00  && 645.900 &57.500 & 703.400 &&3.01387\\
					$(2^{13}, 2^{14})$& 6.19e-05& 1.00 &3.89e-04& 1.00  && -	&-	&-	&&16.6474 \\
					$(2^{14}, 2^{15})$& 3.09e-05& 1.00 &  1.94e-04& 1.00  && -	&-	&-	&&138.398 \\
					\hline
		\end{tabular}}}}
	\end{table}
	
	\begin{figure}[htb]
		\hspace*{\fill}\begin{minipage}[t]{0.45\textwidth}
			\centering
			\includegraphics[width=\textwidth]{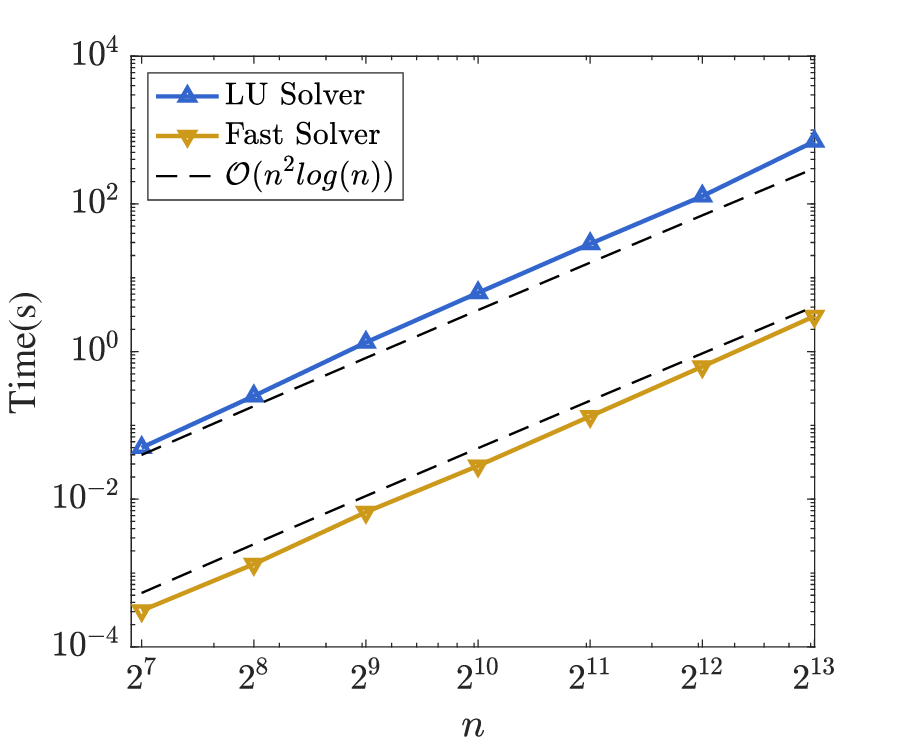}
			\vspace{-2em}
			\caption{\small{Total computational time versus $n$  in Example \ref{exam2}.}}\label{fig_exam2}
		\end{minipage}\hfill%
		\begin{minipage}[t]{0.45\textwidth}
			\centering
			\includegraphics[width= \textwidth]{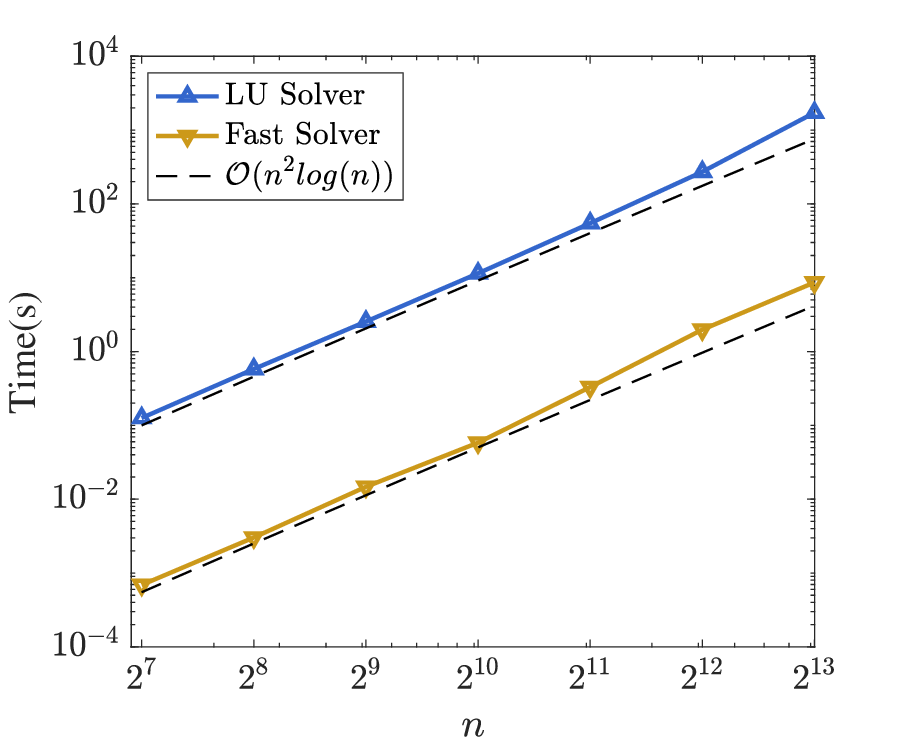}
			\vspace{-2em}
			\caption{\small{Total computational time versus $n$ in Example \ref{exam3}.}}\label{fig_exam3}
		\end{minipage}\hspace*{\fill}
	\end{figure}
	
	Firstly, Table \ref{tab_exam2-1} presents the errors and convergence orders in the norm $\|\bm e_h\|$ and $\|{\color{black} \rot\ }\bm e_h \|$ for these two methods with different partition sizes $n_x$ and $n_y$ in the $x$- and $y$-direction, respectively. When
	$(n_x, n_y)=(2^{12},2^{13})$, both  the {\color{black}LU solver} and  the fast solver
	can reach a correct solution with perfect convergence orders. 
	For a larger partition size, only  the  fast solver  achieves the satisfactory solution.
	The asymptotic behaviors of computational time shown in Figure \ref{fig_exam2} are similar to  those in Example \ref{exam1}. 
	The fast solver requires less computational time than  the  {\color{black}LU solver. 
 The time of the faster solver grows in the order of  $\mathcal{O}(n^2\log n)$, while the time of the LU solver grows faster than  $\mathcal{O}(n^2\log n)$.}
	

	\begin{example}\label{exam3}
		Consider the Maxwell's equation
		\eqref{Maxwell_source}
		when $\alpha=1$ with homogeneous natural boundary conditions \eqref{natural_BCs}. The right-hand side function $\bm{f}$ is chosen such that the exact solution is
		$$\bm{u} =\big(\sin(\pi x)\cos(\pi y), -\cos(\pi x) \sin(\pi y)\big)^{\tr}.$$		
	\end{example}

	\begin{table}[htb]
		\caption{\small{Errors, convergence rates and computational time of {\color{black}LU solver} and fast solver in Example \ref{exam3}.  }}\label{tab_exam3-1}
		\tabcolsep =0.35em
		\vspace{0.1cm}
		\centerline{\small{
				\begin{tabular}{c|ccccc||cccc|c}
					\hline
					\quad & \multicolumn{4}{c}{ Errors and convergence order} & &\multicolumn{5}{c} {Computational time}\\
					\hline
					\quad & \multicolumn{4}{c}{ Fast Solver} & &\multicolumn{3}{c} {{\color{black}LU Solver}}&& Fast Solver\\
					\hline
					$n$ & $\|\bm e_h \|$ & order & $\|{\color{black} \rot\ }\bm e_h \|$ & order& & {\color{black}fact.}(s) & {\color{black}solution}(s)& total(s) && total(s)\\
					\hline
					$2^7$ &5.01e-03&  &3.15e-02&  && 0.120 & 0.006 & 0.126 & &0.00069\\
					$2^8$ &2.50e-03& 1.00 &1.57e-02& 1.00 &&0.546 & 0.034 & 0.580 && 0.00303\\
					$2^9$ &1.25e-03& 1.00 &7.87e-03& 1.00  &&2.370 & 0.178 & 2.548 && 0.01472\\
					$2^{10}$&  6.26e-04& 1.00 & 3.93e-03&1.00  &&  10.776 & 0.684 & 11.460 && 0.05832 \\
					$2^{11}$ & 3.13e-04	&1.00 & 1.97e-03& 1.00  && 51.718 &  2.990 & 54.708 && 0.33000\\
					$2^{12}$& 1.57e-04	& 1.00 & 9.84e-04& 1.00  && 259.160 & 13.000 & 272.160 && 1.98000  \\
					$2^{13}$& 7.83e-05	&1.00 & 4.92e-04& 1.00  && 1615.560 & 103.970 & 1719.530 && 8.63000\\
					$2^{14}$& 3.91e-05& 1.00 & 2.46e-04& 1.00  && -	&-	&-	&&41.2900  \\
					$2^{15}$& 1.96e-05& 1.00 & 1.23e-04& 1.00  && -	&-	&-	&&207.580 \\
					\hline
		\end{tabular}}}
	\end{table}
	
	In such a case,
	similar conclusions could be  drawn from Table \ref{tab_exam3-1}.
	For both  the {\color{black}LU solver} and the fast solver, the errors  in  the $\bm L^2$-norm and {\color{black}$H(\tc;\Omega)$-semi norm} are almost  the same for $n$ ranging from $2^7$ to $2^{13}$, and the corresponding convergence orders agree with the theoretical orders. 
	 However, when $n\ge 2^{14}$, only  the fast solver obtains a valid solution.
	 Regarding computational time, Figure \ref{fig_exam3} indicates that the fast solver  is more efficient than the {\color{black}LU solver}.
	In fact,   Table \ref{tab_exam3-1}  also reveals that the total time used by the {\color{black}LU solver} is 1719.53 seconds, compared  to 8.63 seconds for  the fast solver when $n=2^{13}.$

	\begin{example}\label{exam4}
		Consider the Maxwell's equation with variable coefficients and homogeneous essential boundary conditions
		\begin{align}
			\begin{cases}
				{\color{black} \mathrm{\curl}\ \left(\beta\  \rot\ \bm{u} \right)} + \alpha\bm{u}= \bm{f},\quad &{\rm in}\; \Omega,\\
				\bm{u} \times \bm{n} = 0, \quad &{\rm on}\; \partial\Omega,
			\end{cases}
		\end{align}
  {\color{black}where $\beta\ge \beta_0, \alpha>0$}. 
		The right-hand side function $\bm{f}$ is chosen such that the exact solution is
		$$\bm{u} = \big(\sin^2\pi x\sin \pi y\cos \pi y,-\sin^2\pi y\sin \pi x\cos \pi x\big)^{\tr},$$ while $\beta$ and $\alpha$ are functions with
		\begin{align*}
			&\beta=3\pi\cos(\pi x)\cos(\pi y)+10,
			\quad\alpha=3\pi\sin(\pi x)\sin(\pi y).
		\end{align*}
	\end{example}
	Since the divergence-free condition is guaranteed by the  specifically chosen $\bm f$  when $\alpha>0$ in $\Omega=(0,1)^2$, the Lagrange multiplier is not necessary, resulting in a linear algebraic system  that is symmetric and positive definite.  
	We   will initially employ CG method to solve the  linear equations and then  examine the behaviors of PCG method. The discrete matrix calculated when $\beta=1, \alpha=1$ is chosen as the pre-conditioner, enabling the use of FSTs/FCTs in matrix-vector multiplications. The {\color{black}LU solver}  is also produced  for comparison with PCG.
	
\begin{figure}[htb]
		\begin{minipage}[t]{0.45\linewidth}
			\hfill%
			\includegraphics[width=\textwidth]{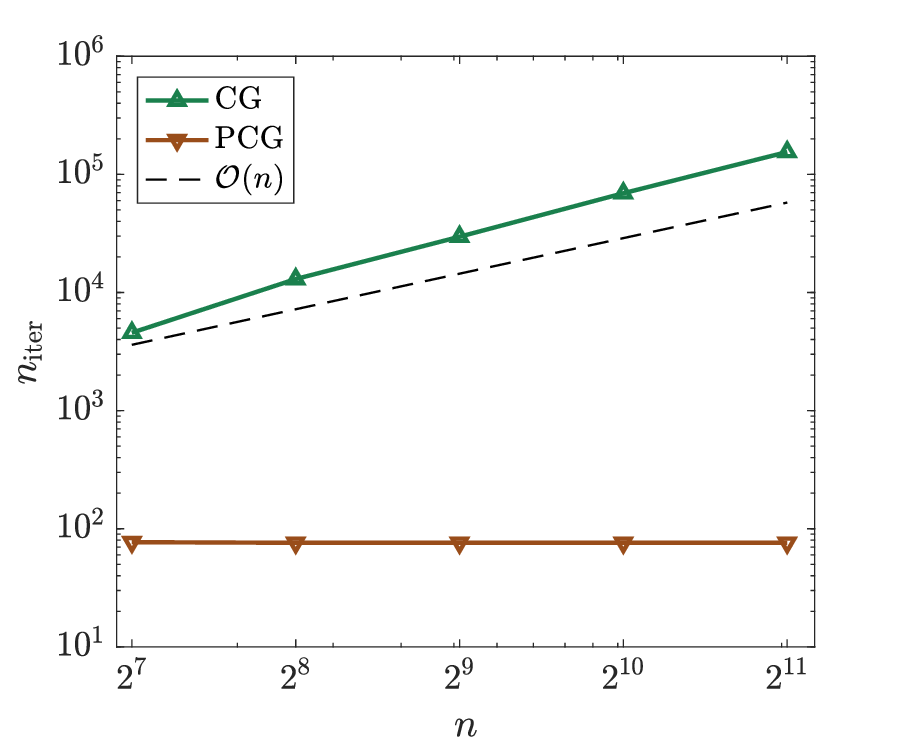}
			\vspace{-2em}
			\caption{\small{Iterations versus $n$ for CG and PCG  in Example \ref{exam4} with $\varepsilon=10^{-14}$.}}\label{fig_exam4-1}
		\end{minipage} \hfill%
		\begin{minipage}[t]{0.45\linewidth}
			\centering
			\includegraphics[width=\textwidth]{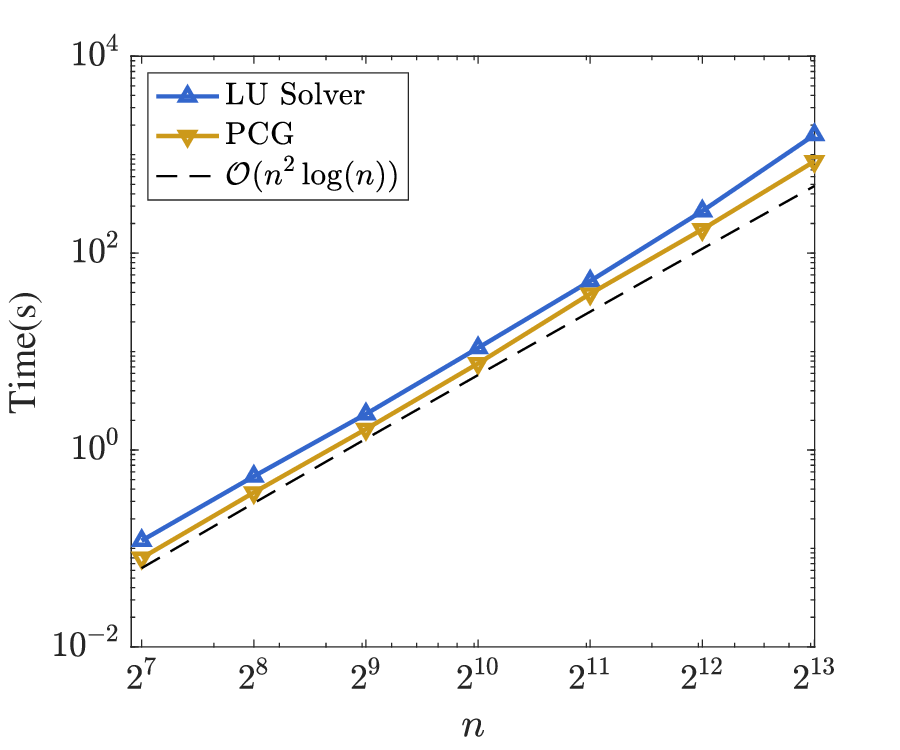}
			\vspace{-2em}
			\caption{\small{Total computational time versus $n$ by PCG  in Example \ref{exam4} with $\varepsilon=10^{-14}$.}}\label{fig_exam4-2}
		\end{minipage}\hspace*{\fill}		
	\end{figure}

	Denote the number   of iterations by $n_{\rm iter}$.
	The exit criterion for the iteration relies  on reducing the 2-norm of the relative residual below the prescribed tolerance $\varepsilon=10^{-14}$.
	In Table \ref{exam4_1e-12}, it is demonstrated that errors in  the $\bm L^2$-norm obtained by the three methods (CG,  PCG, and {\color{black}LU solver}) converge at the rate of $\mathcal{O}(h)$,  consistent with the theoretical results.
	 As $n$ increases, the number  of iterations for the CG method grows faster than the rate of $\mathcal{O}(n)$
	while that for the PCG method is independent of mesh sizes. 
	This   is explicitly illustrated in Figure \ref{fig_exam4-1},
	 confirming the effectiveness of the preconditioner.
	 Additionally, the total time of the  three  methods is recorded  in Table \ref{exam4_1e-12}. 
	For PCG and the {\color{black}LU solver}, we  also record factorization phase time and the preconditioned phase time, respectively. 
	 In contrast, the computational time of PCG  shows orders of magnitude improvement.
	Moreover,  regarding computational time, we find from Figure \ref{fig_exam4-2}
	that for  the PCG method, it grows asymptotically as $\mathcal{O}(n^2 \log n).$ Thus, besides  the reduction  in the number of  iterations, the usage of FSTs/FCTs can also account for the saved computational time.
	 Memory usage  for PCG and  the {\color{black}LU solver}  is also recorded  in Table \ref{mem_exam4-1}. When $n=2^{13},$  the {\color{black}LU solver} requires 57.4\% of total memory, while the fast solver only needs 9.2\%.
{\color{black} The PCG method performs matrix iterations, since the coefficient matrix needs to be stored, which requires additional memory. In contrast, our fast solver is  matrix-free, so PCG uses more memory than the fast solver used in previous examples.
}
	\begin{table}[htb]
		\caption{\small{Errors, convergence rates, iterations and computational time of CG, PCG and {\color{black}LU solver} in Example \ref{exam4} with  $\varepsilon=10^{-14}$. }}\label{exam4_1e-12}
		\tabcolsep =0.3em
		\vspace{0.1cm}
		\centerline{\small{
				\begin{tabular}{c|ccc||ccc||cccc||cccc}
					\hline
					\quad& \multicolumn{2}{c}{Errors \& order}&&\multicolumn{2}{c}{CG} && \multicolumn{3}{c}{PCG} && \multicolumn{3}{c}{{\color{black}LU Solver}} \\
					\hline
					$n$& $\|\bm e_h \|$ & order && $n_{\rm iter}$ & time(s)& &$n_{\rm iter}$ & {\color{black}precond.}(s)& total(s)&&  {\color{black}factor.}(s) & total(s) \\
					\hline
					$2^7$ & 2.51e-3	&&& 4564 & 1.64 && 77&0.05&	0.08&&0.110&0.120\\ 
					$2^8$ & 1.25e-3	& 1.01& & 12994 & 21.19 && 76	&0.24& 0.37&&0.510&0.540\\ 
					$2^9$ & 6.26e-4	& 1.00& & 29654& 202.11 & & 76&	1.10& 1.63&&2.170&2.320\\ 
					$2^{10}$ & 3.13e-4	& 1.00& & 69235 & 2363.02 &&76&4.89& 7.57&&10.200&10.880\\ 
					$2^{11}$ & 1.57e-4& 1.00&& 154719& 26587.23&& 76&26.29&38.62 &&49.170&52.040\\
					$2^{12}$ &7.88e-5& 0.99&&  -&   -&& 76& 124.28 &   174.12 &&254.370&267.480\\ 
					$2^{13}$ &3.99e-5& 0.98&& - &   - && 76  & 651.54& 856.05&&1508.950&1591.490&\\ 
					$2^{14}$ &2.11e-5 & 0.92&& - &   - && 76  &3277.46& 4094.05&&- &-&\\ 
					\hline
		\end{tabular}}}
	\end{table}

	\begin{table}[htb]
		\caption{\small{Memory of {\color{black}LU solver} and PCG in Example \ref{exam4}. (total memory: 256GB) }}\label{mem_exam4-1}
		\centerline{\small{
				\begin{tabular}{c|ccccccccc}
					\hline
					$n$  &$2^9$  & $2^{10}$ & $2^{11}$ &  $2^{12}$ &$2^{13}$&$2^{14}$  \\
					\hline
					mem(PCG) (MB)&0.0&262 & 1311& 6029&24117&95945\\
					{\color{black}mem(LU Solver)} (MB) &262& 1573&8389&35652& 150470&- \\
     {\color{black}mem(PCG)/mem(LU Solver) (\%)}  &0.0& 16.66&15.63&16.91& 16.03&- \\
					\hline
		\end{tabular}}}
	\end{table}
	
	\begin{example}\label{exam5}
		Consider the Maxwell's equation
		\eqref{Maxwell_source_div}-\eqref{Maxwell_source_boundary}
		when $\alpha=1$. The right-hand side function $\bm{f}$ is chosen such that the exact solutions are
		$$\bm{u} = \big(\cos(\pi x)\sin(\pi y), \sin(\pi x) \sin(\pi y)\big)^{\tr},\quad p = \sin(\pi x)\sin(\pi y).$$		
	\end{example}
	
	To  facilitate comparison,  a {\color{black}LU solver} is also performed in this case.
	In Table \ref{tab_exam5-1}, we  observe  that the errors of $\bm e_h$ in  the $\bm L^2$-norm and {\color{black}$H(\tc;\Omega)$-semi norm} for these two methods are almost  the same,  with the convergence order being $\mathcal{O}(h)$,
	 consistent with the theoretical results.
	 Despite both methods exhibiting asymptotic behaviors of total time growing in the same order $O(n^2 \log n)$, as illustrated in Figure \ref{tab_exam5-1}, the fast solver significantly reduces the computational time for the same linear algebraic system.
	 This  once again verifies the efficiency of the fast solver.

	\begin{table}[htb]
		\caption{\small{Errors, convergence rates and computational time of {{\color{black}LU solver}} and  fast solver in Example \ref{exam5}}.  }\label{tab_exam5-1}
		\tabcolsep =0.3em
		\vspace{0.1cm}
		\centerline{\small{
				\begin{tabular}{c|ccccc||cccc|c}
					\hline
					\quad & \multicolumn{4}{c}{ Errors and convergence order} & &\multicolumn{5}{c} {Computational time}\\
					\hline
					\quad & \multicolumn{4}{c}{ Fast Solver} & &\multicolumn{3}{c} {{\color{black}LU Solver}}&& Fast Solver\\
					\hline
					$n$ & $\|\bm e_h \|$ & order & $\|{\color{black} \rot\ }\bm e_h \|$ & order& &  {\color{black}factor.}(s) & {\color{black}solution}(s)& total(s) && total(s)\\
					\hline
					$2^7$ &5.01e-03& &2.23e-02& &&0.300 & 0.026 & 0.326 && 0.0012 \\
					$2^8$ &2.50e-03& 1.00 & 1.11e-02	& 1.01 &&1.666 & 0.158 & 1.824 	&&0.0053 \\
					$2^9$ &1.25e-03& 1.00 & 5.56e-03& 1.00 && 9.246 & 0.772 & 10.018 	&&0.0240 \\
					$2^{10}$&  6.26e-04& 1.00 & 2.78e-03	&1.00 && 44.786 & 3.774 &	48.560 &&0.1080  \\
					$2^{11}$ & 3.13e-04 & 1.00 & 1.39e-03 & 1.00 && 230.455 & 13.995 	& 244.450 &&0.5420 \\
					$2^{12}$& 1.57e-04& 1.00 & 6.96e-04& 1.00 && 1310.510 & 67.160 	& 1377.670 &&2.5180  \\
					$2^{13}$& 7.83e-05&1.00 & 3.48e-04& 1.00 &&-& -& -&&12.4820 \\
					$2^{14}$& 3.91e-05	& 1.00 & 1.74e-04& 1.00 	&&-	&-	&-	&&60.3300  \\
					$2^{15}$&1.96e-05	& 1.00 & 8.69e-05& 1.00 &&	-& -	&-&&	266.2700  \\
					\hline
		\end{tabular}}}
	\end{table}

	\begin{figure}[htb]
		\begin{center}
			\includegraphics[width=6cm]{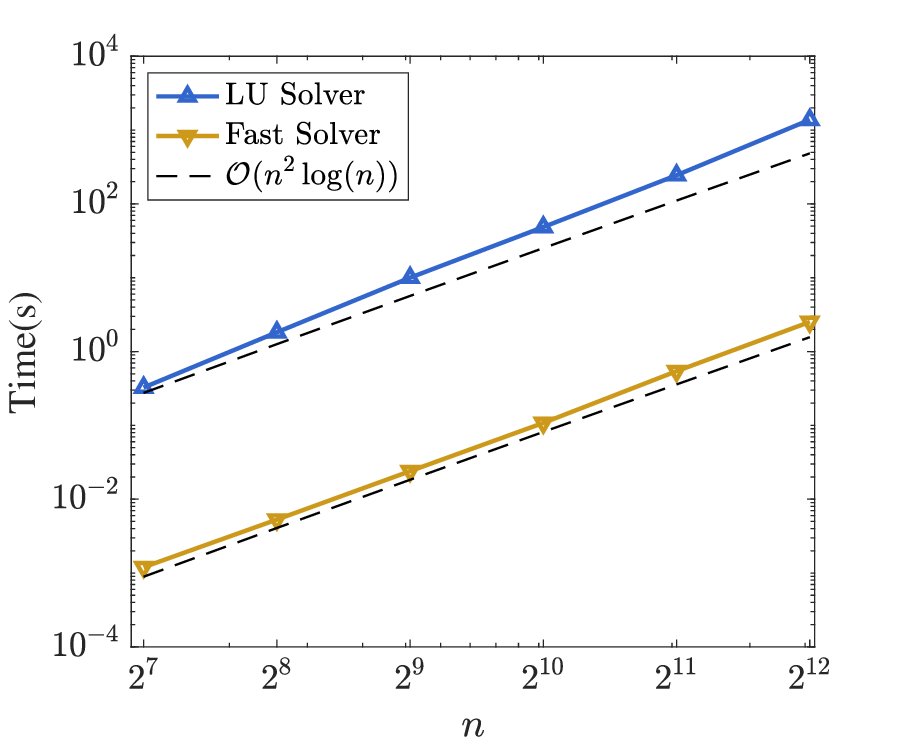}
			\vspace{-1em}
			\caption{Total computational time versus $n$ by Fast Solver in Example \ref{exam5}.}\label{fig_exam5-2}
			\label{convergenceR}
		\end{center}
	\end{figure}

	\section{Conclusion}\label{conclusion}
	In this paper,   fast solvers for the Maxwell's equations are  proposed by using the divergence-free constraint / Gauss's law guaranteed  discrete eigen-decompositon  and the FSTs/FCTs in $\mathcal{O}( n^2\log n)$ operations when adopting 
	the lowest-order N\'ed\'elec finite element.  Notably, these algorithms strictly fulfill the Helmholtz-Hodge decomposition at the discrete level, thus instabilities and  spurious solutions can be  eliminated.
  {\color{black}Moreover, based on ideas from Isogeometric Analysis, we can extend our
fast solver to more general computational domains $\Omega$, which  can be easily parameterized by a sufficiently smooth function
$G : (0,1)^2 \rightarrow \Omega$. In this case, the variational problem on $\Omega$ could be transferred to $(0,1)^2$
using substitution and chain rules for integration and differentiation.}  
			 Additionally, as part of an ongoing series of works, we are considering the design of fast solvers for Maxwell's equations in  three-dimensions, which  offer substantial advantages over existing methods.



\section*{Acknowledgements}
This work is supported in part  by  the National Key R\&D Program of China (No. 2021YFB0300203), 
the National Natural Science Foundation of China grants (NSFC  12101036, NSFC   12101035, NSFC 12471348, and NSFC 12131005),
and the Fundamental Research Funds for the Central Universities (No. FRF-TP-22-096A1).
\\
\bigskip

{\noindent \bf {Declarations of interest: none}}

\bibliographystyle{elsarticle-num-names}

\bibliography{reference}

\end{document}